\documentclass[12pt]{elsarticle}
\usepackage{amssymb, amsmath, amsfonts}
\usepackage{hyperref}
\usepackage[all]{xy}
\usepackage{bm}
\def\genfd{{\bm k}}
\def\op{\mathrm{op}}
\def\id{\mathrm{id}}
\long\def\nodo#1{{}}
\def\gg{\mathfrak{g}}

\def\hx{\hat{x}}

\def\MR#1{} 
\def\hx{\hat{x}}

\def\btr{\blacktriangleright}

\usepackage{amsthm}
\newtheoremstyle{definition}{}{}{\upshape}{}{\bfseries}{.}{0.5em}{}
\theoremstyle{definition}
\newtheorem{theorem}{Theorem}[section]
\newtheorem{lemma}[theorem]{Lemma}
\newtheorem{corollary}[theorem]{Corollary}
\newtheorem{proposition}[theorem]{Proposition}
\newtheorem{definition}[theorem]{Definition}

\newtheorem{remark}[theorem]{Remark} 
\newtheorem{notation}[theorem]{Notation}

\begin{document}
\title{Hopf algebroids with balancing subalgebra}
\author[uhk]{Zoran \v{S}koda\corref{cor}}\ead{zoran.skoda@uhk.cz}
\author[pmfmo]{Martina Stoji\'c}\ead{stojic@math.hr}
\cortext[cor]{Corresponding author}
\address[uhk]{Faculty of Science, University of Hradec Kr\'{a}lov\'{e}, Rokitansk\'{e}ho 62, Hradec Kr\'{a}lov\'{e}, Czech Republic}
\address[pmfmo]{Department of Mathematics, University of Zagreb,
Bijeni\v{c}ka cesta~3, HR-10000 Zagreb, Croatia}
\begin{abstract}
Recently, S.\ Meljanac proposed a construction of a class
of examples of an algebraic structure with properties very close
to the Hopf algebroids $H$ over a noncommutative base $A$ of other authors.
His examples come along with a subalgebra $\mathcal{B}$ of $H\otimes H$,
here called the balancing subalgebra, which  
contains the image of the coproduct and such that 
the intersection of $\mathcal{B}$ with the 
kernel of the projection $H\otimes H\to H\otimes_A H$ 
is a two-sided ideal in $\mathcal{B}$ 
which is moreover well behaved with respect to the antipode. 
We propose a set of abstract axioms covering this construction
and make a detailed comparison to the Hopf algebroids of Lu.
We prove that every scalar extension Hopf algebroid can be
cast into this new set of axioms. We present an observation
by G.\ B\"ohm that the Hopf algebroids constructed
from weak Hopf algebras fit into our framework as well.
At the end we discuss the change of balancing subalgebra under
Drinfeld-Xu procedure of twisting of associative bialgebroids
by invertible 2-cocycles.
\end{abstract}
  \maketitle

\section{Introduction}

Hopf algebroids~\cite{bohmHbk,BrzMilitaru,lu} 
are generalizations of Hopf algebras,  
which are roughly in the same relation to groupoids as Hopf algebras 
are to groups. 
They are {\bf bialgebroids} possessing a version of an antipode,
where an (associative) bialgebroid is 
the appropriate generalization of a bialgebra.
Hopf algebroids comprise several structure maps 
defined on a pair of associative unital algebras, 
the {\bf total algebra} $H$ (generalization of a function algebra on 
the space of morphisms of a groupoid), 
and the {\bf base algebra} $A$ (generalization of a function algebra 
on the space of objects (equivalently, units) of a groupoid). 
The main structure on the total algebra of a bialgebroid is an $A$-bimodule
structure on $H$ and a coproduct $\Delta:H\to H\otimes_A H$.
The commutative Hopf algebroids (where both $H$ and $A$ are 
commutative) are easy to define by a categorical dualization of
the groupoid concept.
They are used as a classical tool in 
stable homotopy theory \cite{hovey,Ravenel}.
Noncommutative Hopf algebroids over
a commutative base ($H$ noncommutative, $A$ commutative) are also rather 
straightforward to introduce; this theory has been studied from late 1980-s, 
under the influence of the quantum group theory~\cite{maltsiniotis}. 
The most obvious examples are the convolution algebras of finite groupoids. 
Bialgebroids and Hopf algebroids 
over a noncommutative base are much more complicated to
define; several versions were developed in early 1990s by 
Lu~\cite{lu}, Xu~\cite{xu},
B\"ohm~\cite{bohmnew}, B\"ohm-Szlach\'anyi~\cite{BohmSzlach},
Day-Street (in~\cite{daystreetHopfalgd} and a different abstract notion in~\cite{daystreetqcat}) and others,
including an earlier 
notion of $\times_A$-bialgebra~\cite{takeuchi} and its Hopf version by Schauenburg~\cite{schauenburg};
for comparisons see~\cite{bohmHbk,BrzMilitaru}. B\"ohm has also introduced
an internalization of a bialgebroid in any symmetric monoidal category with coequalizers commuting with tensor product~\cite{bohmInternal}; this has been
extended to an internalization of Hopf algebroids in~\cite{stojicPhD}.
Many examples of Hopf algebroids over
noncommutative bases have been studied in the contexts of 
inclusions of von Neumann algebra factors~\cite{kadszl,bohmHbk},
dynamical Yang-Baxter equation~\cite{doninmudrov}, weak Hopf algebras~\cite{bohmHbk}, 
deformation quantization~\cite{xu},
noncommutative torsors, noncommutative differential calculus and
cyclic homology~\cite{kow} etc. 

In 2012, S.\ Meljanac devised a new approach to some examples of
(topological) Hopf algebroids over a noncommutative base 
restricting the codomain of the coproduct 
map in a useful, but somewhat {\it ad hoc} way. To construct that
codomain, he chooses a subalgebra $\mathcal{B}$ in the tensor
square $H\otimes H$ of the total algebra $H$,
such that the intersection of $\mathcal{B}$ with the kernel $I_A$ of
the projection $\pi: H\otimes H\to H\otimes_A H$ 
to the tensor square over the noncommutative base algebra $A$
is a two-sided ideal $I_A\cap\mathcal{B}\subset\mathcal{B}$
(with an appropriate behaviour under the antipode map). 
The appearance of the two-sided ideal is a novel
and somewhat unexpected feature reminding of the classical case where
the base algebra $A$ is commutative and $I_A$ is two-sided itself.
The approach is developed in collaborative works~\cite{tajron},  
with more details in~\cite{tajronkov}. 

Papers~\cite{tajronkov} and~\cite{tajron}
neglect two mathematical issues. 
Firstly, no care is taken about the implicit use of completions: 
the values of the coproduct involve infinite sums, 
hence its codomain should be a completed tensor square. 
Secondly, at the algebraic level, they do not state a complete
axiomatic framework for their version of Hopf algebroid, nor state
its precise relation to other definitions.  
Instead, they construct an interesting class of examples
and give a partial list of essential properties. 
In our article~\cite{halgoid} with Meljanac, using the explicit formulas
from~\cite{ldWeyl}, we treat 
a somewhat wider class of examples 
in a mathematically rigorous way, 
using G.~B\"ohm's definition of a symmetric Hopf algebroid, 
{\em partly} adapted to a formally completed tensor product. 
For a better adaptation, which gives rise to an internal Hopf algebroid
in a symmetric monoidal category of filtered-cofiltered vector spaces,
entailing a more sensible completion, see~\cite{stojicPhD,stojicindproVect}. These works took care of completions,
but instead of the two-sided ideal approach they relied on (an internalization of) symmetric Hopf algebroid axiomatics~\cite{bohmHbk,BohmSzlach}. 
To return closer to the original idea,
we here propose a new set of axioms 
expressing the essence of the two-sided ideal approach and discuss
it in the context. The subalgebra $\mathcal{B}\subset H\otimes H$ 
in new axioms is named the {\bf balancing subalgebra} and our new
version of Hopf algebroid over noncommutative base algebra $A$
is named a Hopf $A$-algebroid with balancing subalgebra $\mathcal{B}$.

In Theorem~\ref{th:MLu} we compare Hopf A-algebroid with
a balancing subalgebra to the Hopf algebroids of Lu instead
to symmetric Hopf algebroids from~\cite{bohmHbk,BohmSzlach,halgoid}.
This is because Lu's axioms for the antipode map 
involve a choice of certain map (section $\gamma$ below) which is 
close in spirit to the choice of balancing subalgebra
in our axiomatics and in the informal approach of Meljanac.
Our main result is Theorem~\ref{th:scalarM}
(based on nontrivial Lemmas~\ref{lem:key1},\ \ref{lem:key2})
stating that every scalar extension Hopf algebroid 
can be cast into a Hopf algebroid
with a suitable choice of balancing subalgebra $\mathcal{B}$.

In Section~\ref{sec:wha}, we worked out the observation of B\"ohm that
each weak Hopf algebra gives rise not only to Hopf algebroids in the senses
of Lu~\cite{lu} and B\"ohm~\cite{bohmnew,BohmSzlach},
but also to a Hopf algebroid with a balancing subalgebra.

Throughout the paper, $\genfd$ is a commutative ground field, and
the unadorned tensor symbol $\otimes$ between symbols
for $\genfd$-vector spaces, is meant over $\genfd$,
however we often use $\otimes_\genfd$ for emphasis.
When used among elements in calculations,
symbol $\otimes$ is interpreted from the context.

\section{Bialgebroids over noncommutative base}

The axioms of bialgebroids and Hopf algebroids over a noncommutative base algebra are far less obvious to formulate~\cite{bohmHbk,BrzMilitaru,lu,schauenburg,xu}. Let us detect a problem naively.
For a commutative $\genfd$-algebra $A$,
an $A$-bialgebra is a monoid (algebra) 
and a comonoid (coalgebra) in the same symmetric 
monoidal category, namely that of $A$-modules, with a compatibility condition utilizing the symmetry of the tensor product $\otimes_A$.
For a noncommutative base algebra $A$ over $\genfd$,
the category of $A$-modules is not monoidal,
so it is natural to try replacing it
with the monoidal category of $A$-bimodules. 
However, the latter is neither symmetric nor braided monoidal
in general, so the usual compatibility condition 
between the comonoids and monoids makes no sense. 
Instead, it appears that 
the monoid and the comonoid part of the left $A$-bialgebroid
structure live in different monoidal categories~\cite{bohmHbk}. 

The monoid structure $(H,\mu,\eta)$ on $H$ is in 
the monoidal category of $A \otimes A^\op$-bimodules;
equivalently (\cite{bohmHbk}, Lemma 2.2),
$(H,\mu)$ is an associative $\genfd$-algebra 
and the unit map $\eta$ is a morphism of $\genfd$-algebras 
$\eta: A\otimes_\genfd A^\op\to H$ 
(we say that $(H,\mu,\eta)$ is an $A\otimes A^\op$-ring).
The unit $\eta:A\otimes A^\op\to H$ is usually described in terms of its 
left leg $\alpha := \eta(-\otimes 1_{A^\op}):A\to H$ and
its right leg $\beta := \eta(1_A\otimes -):A^\op\to H$, 
also called the source and target maps respectively; then,
their images commute because
\begin{equation}\label{eq:alphabetacom}
\alpha(a)\beta(b) = 
\eta(a\otimes 1)\eta(1\otimes b) = \eta(a\otimes b)=
\eta(1\otimes b)\eta(a\otimes 1) =
\beta(b)\alpha(a).
\end{equation}
An $A\otimes A^\op$-ring $(H,\mu,\eta)$ is 
described below as the equivalent datum $(H,\mu,\alpha,\beta)$.

On the other hand, the comonoid structure $(H,\Delta,\epsilon)$ 
is in the monoidal category
of $A$-bimodules (we say that $H$ is an $A$-coring,~\cite{BrzWis}). 

\begin{definition} 
An $A\otimes A^\op$-ring $(H,\mu,\alpha,\beta)$ and 
an $A$-coring $(H,\Delta,\epsilon)$
with underlying $A$-bimodule $H$
form a {\bf left associative $A$-bialgebroid} 
$(H,\mu,\alpha,\beta,\Delta,\epsilon)$
if they satisfy the following compatibility conditions: 
\begin{itemize}
\item[(C1)] The underlying $A$-bimodule
structure of $A$-coring $H$ is determined by the source and
target maps (part of the $A\otimes A^\op$-ring structure): 
$a.h.b = \alpha(a)\beta(b)h$ for $a,b\in A$ and $h\in H$.
This is indeed a bimodule by~(\ref{eq:alphabetacom}).
\item[(C2)] Formula
  $\triangleright\colon\sum_\lambda h_\lambda\otimes a_\lambda \mapsto
\epsilon(\sum_\lambda h_\lambda\alpha(a_\lambda))$ defines an action 
$H\otimes A\stackrel{\triangleright}\to A$.
\item[(C3)] The map
  $H\otimes_\genfd (H\otimes_\genfd H)\to H\otimes_A H$ given by the rule
  $h\otimes (g\otimes f)\mapsto \Delta(h)(g\otimes f)$,
factorizes through a map $H\otimes_\genfd (H\otimes_A H)\to H\otimes_A H$
which is moreover a unital action.
Expression $\Delta(h)(g\otimes_\genfd f)$ is understood by taking any representative of $\Delta(h)$ in $H\otimes_\genfd H$, then multiplying in each tensor factor separately with $g\otimes f\in H\otimes_\genfd H$; the result is well defined in
$H\otimes_A H$. Unitality of the action implies $\Delta(1) = 1\otimes_A 1$.
\end{itemize}
\end{definition}

By (C1), $\epsilon$ being a bimodule map means $\epsilon(\alpha(a)h) = a\epsilon(h)$ and $\epsilon(\beta(b)h) = \epsilon(h)b$. In particular, $\epsilon\circ\alpha = \epsilon\circ\beta = \id_A$ and $\epsilon(1_H) = \epsilon(\alpha(1_A)) = 1_A$. Action in (C2) is unital by its defining formula and it extends the left regular action $A\otimes A\to A$ along the inclusion 
$A\otimes A\stackrel{\alpha\otimes A}\longrightarrow H\otimes A$ by direct calculation, $\epsilon(\alpha(a)\alpha(b)) = \epsilon(\alpha(a b)) = a b$.
Action axiom (C2), for acting on $1_A$,
implies $\epsilon(h g) = \epsilon(h g\alpha(1_A)) = (h g)\triangleright 1_A =
h\triangleright (g\triangleright 1_A) = \epsilon(h\alpha(\epsilon(g)))$. In
particular, $\epsilon(h\beta(b)) = \epsilon(h\alpha((\epsilon\circ\beta)(b)))
= \epsilon(h\alpha(b)) = h\triangleright b$. Action axiom on $1_A$ together with $a = \epsilon(\alpha(a)) = \alpha(a)\triangleright 1$, implies the general case, $h\triangleright(g\triangleright a) = h\triangleright (g\triangleright(\alpha(a)\triangleright 1)) = h\triangleright((g\alpha(a))\triangleright 1) =
(h g \alpha(a))\triangleright 1 = (h g)\triangleright(\alpha(a)\triangleright 1)
= (h g)\triangleright a$.

From (C1) and $\Delta(1_H) = 1_H\otimes_A 1_H$, $\Delta$ being an $A$-bimodule map implies $\Delta(\alpha(a)) = \Delta(\alpha(a) 1_A) = \alpha(a)\otimes_A 1_H$ and $\Delta(\beta(b)) = 1_H\otimes_A\beta(b)$. It follows that $\Delta(h\alpha(a)) = \Delta(h)(\alpha(a)\otimes_A 1) = h_{(1)}\alpha(a)\otimes_A h_{(2)}$ and $\Delta(h\beta(a)) = h_{(1)}\otimes_A h_{(2)}\beta(a)$. Applying the counit axiom we obtain, for all $h\in H$ and $a\in A$,
\begin{equation}\label{eq:halpha}
  h\alpha(a) = \alpha(\epsilon(h_{(1)}\alpha(a)))h_{(2)} = \alpha(h_{(1)}\triangleright a)h_{(2)},
\end{equation}
\begin{equation}\label{eq:hbeta}
  h\beta(a) = \beta(\epsilon(h_{(2)}\beta(a)))h_{(2)}
  = \beta(h_{(2)}\triangleright a)h_{(1)}.
\end{equation}  

The condition (C1) implies that the kernel $I_A = \operatorname{Ker}\pi$
of the projection map
$$\pi: H\otimes_\genfd H\to H\otimes_A H$$ 
of $H$-bimodules is a {\bf right ideal}
in the algebra $H\otimes_\genfd H$, generated by the set of elements
of the form $\beta(a)\otimes 1 - 1\otimes\alpha(a)$:
\begin{equation}\label{eq:IA}
I_A = \left\{\,\beta(a)\otimes_\genfd 1 - 1\otimes_\genfd\alpha(a)\,|\,\,a\in A\,\right\}\cdot(H\otimes_\genfd H)
\end{equation}
Regarding that $I_A$ is a right ideal, and $\Delta(h)$ is defined up to $I_A$, the map  $H\otimes_\genfd (H\otimes_\genfd H)\to H\otimes_A H$ in (C3) is well defined. Its factorization through a map $H\otimes_\genfd (H\otimes_A H)\to H\otimes_A H$ exists iff for every $h$, $\Delta(h)I_A \subset I_A$, which is clearly equivalent to $\Delta(h)(\beta(a)\otimes 1 - 1\otimes\alpha(a)) \in I_A$ for all $a\in A$. Hence $\Delta(h)$ must belong to a set
$$
H\times_A H = \Big\{\sum_i b_i\otimes b'_i \in H\otimes_A H \ \big|\,
\sum_i b_i\otimes b'_i \alpha(a) = \sum_i b_i \beta(a)\otimes
b'_i, \,\,\forall a\in A \Big\},
$$ 
which is a $A$-subbimodule of $H\otimes_A H$~(\cite{sweedler,takeuchi}),
called the Takeuchi product~\cite{BrzMilitaru,bohmHbk}. In these terms, the factorization property from~(C3) is equivalent to the property $\mathrm{Im}\,\Delta\subset H\times_A H$. Another part of~(C3), stating that the induced map is an action, may also be expressed in terms of Takeuchi product. 
A direct check shows that the preimage $\pi^{-1}(H\times_A H) = \{\sum_i b_i\otimes_\genfd b'_i, \sum_i b_i\otimes b'_i\alpha(a)- b_i\beta(a)\otimes b'_i\in I_A\}$ is closed under multiplication; in fact a unital subalgebra. The right ideal $I_A$ is spanned by the elements of the form $\beta(a)d\otimes d' - d\otimes\alpha(a)d'$ and if $\sum_i b_i\otimes b'_i\in H\times_A H$ then
$$
\left(\sum_i b_i\otimes b'_i\right)(\beta(a)d\otimes d' - d\otimes \alpha(a)d')
= \left(\sum_i b_i\beta(a)\otimes b'_i - b_i\otimes b'_i\alpha(a)\right)
(d\otimes d')
$$
and the right hand side clearly belongs to $I_A$. Thus, $I_A\cap\pi^{-1}(H\times_A H)$ is not only a right ideal but a two sided ideal of $\pi^{-1}(H\times_A H)$ showing that $H\times_A H\cong \pi^{-1}(H\times_A H)/(I_A\cap\pi^{-1}(H\times_A H))$ is, unlike $H\otimes_A H$, 
an associative algebra with respect to the componentwise product.
The componentwise rule is not well defined in $H\otimes_A H$  because
it may depend on the chosen representatives in $H\otimes_\genfd H$; 
this is because $I_A$ is only a {\it right} ideal in general.

This discussion shows that (C3) is equivalent to the joint assertion 
of the following two requirements:
\begin{itemize}
\item[(C3a)]  
$\mathrm{Im}\,\Delta\subset H\times_A H$,  
\item[(C3b)]  $\Delta$ as a map
from $H$ to $H\times_A H$ is a homomorphism of algebras. 
\end{itemize}
Of course, (C3b) makes sense only because of (C3a). 
Observe now a commutative diagram of $A$-bimodules:
\begin{equation}\label{eq:diagspacestimes}
\xymatrix{
&&\pi^{-1}(H\times_A H)\ar[r]\ar[d]^{\pi|_{\pi^{-1}(H\times_A H)}}& H\otimes_\genfd H\ar[d]^\pi\\
H\ar[rr]^{\Delta\,\,\, }&& 
H\times_A H\ar[r]& H\otimes_A H
}
\end{equation}
All arrows except those into $H\otimes_A H$ are also homomorphisms of algebras.

The equation 
$\sum_i b_i\otimes b'_i \alpha(a) = \sum_i b_i \beta(a)\otimes b'_i$
for elements in $H\otimes_A H$
is demanded in the quotient, hence it holds only up to elements in $I_A$;
if we take the same equation strictly in $H\otimes_\genfd H$ to cut
some subalgebra (actually a left ideal)
$H\tilde\times H\subset H\otimes_\genfd H$, 
then the projection
$\pi|_{H\tilde\times H}$ maps this subalgebra within $H\times_A H$, 
but is not necessarily onto.
In a categorical language, 
$H\times_A H$ is an end (kind of a categorical limit) 
of a coend (kind of a colimit), not the other way around. 
However, Meljanac in his examples
takes some other subalgebra $\mathcal{B}\subset H\otimes_\genfd H$ 
(not a universal construction) first and then passes to the quotient by 
$\pi|_{\mathcal{B}}$ (hence a colimit), with a result which is still 
an algebra (different from $H\times_A H$). 
To achieve this, he needs that 
\begin{itemize}
\item[(C3MI)]  $I_A\cap\mathcal{B}$ 
is a two-sided ideal in $\mathcal{B}$.
\end{itemize}
In addition, he (implicitly) requires
\begin{itemize}
\item[(C3Ma)]  
$\mathrm{Im}\,\Delta\subset \mathcal{B}/(I_A\cap\mathcal{B})$,
\item[(C3Mb)]  $\Delta$ as a map
from $H$ to $\mathcal{B}/(I_A\cap\mathcal{B})$ is a homomorphism of algebras. 
\end{itemize}

\begin{definition} A {\bf left $A$-bialgebroid with  
balancing subalgebra $\mathcal{B}$} comprises an 
$A\otimes_\genfd A^\op$-ring $(H,m,\alpha,\beta)$ and 
an $A$-coring $(H,\Delta,\epsilon)$
with the same underlying $A$-bimodule $H$
and satisfying (C1) and (C2),
and a (not necessarily unital) subalgebra 
$\mathcal{B}\subset H\otimes_\genfd H$ satisfying
(C3MI), (C3Ma) and (C3Mb).
$\mathcal{B}$ is called the {\bf balancing subalgebra}.
\end{definition}

A left $A$-bialgebroid with balancing subalgebra $\mathcal{B}$
is not necessarily 
a left associative $A$-bialgebroid in the standard sense, 
because (C3) does not always hold. 
However, if $\mathcal{B}$ is the preimage $\pi^{-1}(H\times_A H)$ 
of the Takeuchi product under the natural projection 
$\pi$ then (C3) follows. 
Conversely, given a left associative $A$-bialgebroid
$H = (H,m,\alpha,\beta,\Delta,\epsilon)$, we have presented above that
$\pi^{-1}(H\times_A H)\cap I_A$ is a two sided ideal of the
subalgebra $\pi^{-1}(H\times_A H)$. Therefore, (C3) implies that
it is a balancing subalgebra called 
the {\bf trivial balancing subalgebra} of $H\otimes_\genfd H$. It follows
that in a left associative $A$-bialgebroid any subalgebra of $\pi^{-1}(H\times_A H)$ containing $\pi^{-1}(\operatorname{Im}\,\Delta)$ is also balancing.
Therefore, on the level of bialgebroids, balancing algebras are interesting only when either we can not determine whether $\operatorname{Im}\Delta\subset H\times_A H$, or when it does not hold but there is a balancing subalgebra, which is
in the latter case automatically not a subalgebra of Takeuchi product $H\times_A H$.
However, for Hopf algebroids, as we shall see below,
balancing subalgebras provide more
flexible approach to introducing the antipode than using Lu's section,
while it is technically more compact (less structure and axioms)
than B\"ohm's symmetric Hopf algebroids.

Observe a commutative diagram of $A$-bimodules where all arrows except 
those into $H\otimes_A H$ are homomorphisms of algebras:

\begin{equation}\label{eq:diagspaces1}
\xymatrix{
&&\mathcal{B}\ar[r]\ar[d]^{\pi|_{\mathcal{B}}}& H\otimes_\genfd H\ar[d]^\pi\\
H\ar[rr]^{\Delta\,\,\, }&& \mathcal{B}/(I_A\cap\mathcal{B})\ar[r]&
H\otimes_A H
}
\end{equation}

\begin{proposition}\label{prop:equivifBexists} 
Let $(H,\mu,\alpha,\beta,\Delta,\epsilon)$ be the data defining
an $A\otimes A^\op$-ring and $A$-coring satisfying (C1), (C2) and (C3a). 
Suppose there exist a subalgebra 
$\mathcal{B}\subset H\otimes_\genfd H$ such that (C3MI) and (C3Ma) hold. 
Then (C3b) iff (CM3b). In other words, 
these data define a left $A$-bialgebroid with balancing 
subalgebra $\mathcal{B}$ iff 
they (without $\mathcal{B}$) define a left associative $A$-bialgebroid.  
\end{proposition}

\begin{proof} 
This is a rather simple observation: (C3a) and (C3Ma) together 
imply that $\mathrm{Im}\,\Delta\subset\frac{\mathcal{B}}{I_A\cap\mathcal{B}}\cap H\times_A H$ which
has the structure of a subalgebra of 
$\mathcal{B}/(I_A\cap\mathcal{B})$ and also of $H\times_A H$;
the multiplications in $\mathcal{B}/(I_A\cap\mathcal{B})$
and in $H\times_A H$ 
are both defined componentwise, hence they are equal on the intersection. 
If we assume (C3b), then algebra map $\Delta:H\to H\times_A H$ corestricts
to algebra map $H\to\frac{\mathcal{B}}{I_A\cap\mathcal{B}}\cap H\times_A H$
which postcomposed with inclusion of algebras into $\mathcal{B}/(I_A\cap\mathcal{B})$ is again an algebra map,
hence $\Delta:H\to\mathcal{B}/(I_A\cap\mathcal{B})$ is also
an algebra map, hence (C3Mb) holds. Likewise we infer (C3b) from (C3Mb).
\end{proof}

\section{Hopf algebroids: antipode}

\begin{definition}\label{def:LuHopf} A {\bf Hopf $A$-algebroid} 
in the sense of J-L.~Lu \cite{lu} (or a Lu-Hopf algebroid)
is a left associative $A$-bialgebroid 
$(H,\mu,\alpha,\beta,\Delta,\epsilon)$ with an antipode map $\tau:H\to H$,
which is a linear antiautomorphism satisfying
\begin{eqnarray}
\tau\beta = \alpha\label{eq:taubetaalpha} \\
\mu(\id\otimes_\genfd\tau)\gamma\Delta = \alpha\epsilon\label{eq:gam} \\
\mu(\tau\otimes_A\id)\Delta = \beta\epsilon\tau\label{eq:nogam}
\end{eqnarray}
for some linear section $\gamma: H\otimes_A H\to H\otimes H$ 
of the projection $\pi: H\otimes H\to H\otimes_A H$. 
\end{definition}
The reason for introducing $\gamma$ in (\ref{eq:gam}) is the fact that
$\mu(\id\otimes_A\tau)\Delta$ is not a well
defined map because $\mu(\id\otimes_\genfd\tau)(I_A)\neq 0$ in general.
Indeed, $I_A$ is a linear span of the set of all elements of the form
$\beta(a)h\otimes k-h\otimes\alpha(a)k$, where $a\in A$ and $h,k\in H$,
and 
$\mu(\id\otimes\tau)(\beta(a)h\otimes k-h\otimes\alpha(a)k) = 
\beta(a)h\tau(k) - h\tau(k)\tau(\alpha(a))$
which can be nonzero in general. No such problems occur
with~(\ref{eq:nogam}) because
$$\mu(\tau\otimes\id)(\beta(a)h\otimes k-h\otimes\alpha(a)k) = 
\tau(h)\tau(\beta(a))k - \tau(h)\alpha(a)k 
\stackrel{(\ref{eq:taubetaalpha})}=0.$$ 

\begin{definition}
A {\bf Hopf $A$-algebroid with balancing subalgebra} $\mathcal{B}$
is a left $A$-bi\-al\-ge\-broid 
$(H,\mu,\alpha,\beta,\Delta,\epsilon)$ 
with balancing subalgebra $\mathcal{B}$
together with an algebra antihomomorphism $\tau:H\to H$, 
called the {\bf antipode}, such that 
\begin{eqnarray}\label{eq:tauvanish}
\mu(\id\otimes_\genfd\tau)(I_A\cap\mathcal{B})=0\label{eq:tauideal} \\
\tau\beta = \alpha\label{eq:mtaubetaalpha} \\
\mu(\id\otimes_A\tau)\Delta = \alpha\epsilon\label{eq:mr} \\
\mu(\tau\otimes_A\id)\Delta = \beta\epsilon\tau\label{eq:ml}
\end{eqnarray}
\end{definition}
Two equations are the same as in Definition~\ref{def:LuHopf}:
(\ref{eq:mtaubetaalpha}) is identical to (\ref{eq:taubetaalpha})
and (\ref{eq:ml}) to (\ref{eq:nogam}).
Equation (\ref{eq:mr}) now makes sense because of (\ref{eq:tauideal}).
There is no more need for a choice of a section $\gamma$. Choice
of the subalgebra $\mathcal{B}$ which accomplishes the same. 

\begin{remark}
The map $\mu(\id\otimes_\genfd\tau): h\otimes h'\mapsto h\tau(h')$ 
is $\genfd$-linear, but neither a homomorphism nor
an antihomomorphism of algebras. Hence,
it is not suffi\-cient to check~(\ref{eq:tauideal}) on the algebra 
generators of $I_A\cap\mathcal{B}$, and {\em a fortiori},
on its generators as an ideal in $\mathcal{B}$.
This will be the central difficulty in Section~\ref{sec:scalar}.  
\end{remark}

\begin{theorem} \label{th:MLu}
If a Hopf algebroid with a balancing subalgebra satisfies (C3{a}) 
then it admits a (possibly nonunique) structure of Lu-Hopf algebroid.
\end{theorem}

\begin{proof} Choose a vector space splitting of $H\otimes_A H$
into $\mathrm{Im}\,\Delta$ and a linear complement; for $\gamma$ 
take any linear section 
of the projection $\pi:H\otimes_\genfd H\to H\otimes_A H$
such that values $\gamma(p)$ over all $p\in\mathrm{Im}\,\Delta$
are in $\mathcal{B}$  (this can be done by (C3Ma))
and on the linear complement prescribe any linear choice for $\gamma$, 
for instance $0$. 
Condition (C3b) holds by (CM3b) 
and Proposition~\ref{prop:equivifBexists}.
Then $\mu(\id\otimes_\genfd\tau)\gamma\Delta(h)=\mu(\id\otimes_A\tau)\Delta(h)$
as the right hand side is defined by choosing any representative of $\Delta(h)$
in $H\otimes_\genfd H$ and evaluating $\mu(\id\otimes_\genfd\tau)$. 
Thus~(\ref{eq:gam}) holds, and other conditions on the antipode 
become identities. 
\end{proof}

B\"ohm and Szlach\'anyi exhibited Example\ 4.9 in~\cite{BohmSzlach} of a symmetric Hopf algebroid which does not carry a structure of a Lu-Hopf algebroid. An application of Theorem~\ref{th:MLu} implies that this example is not a Hopf algebroid with a balancing subalgebra either.
We do not know if for any Lu-Hopf algebroid there is a balancing subalgebra (containing the image of $\gamma$). However, we exhibit recipes for a balancing subalgebra for several most prominent classes of Lu-Hopf algebroids. Notably, in Section~\ref{sec:scalar}, for any scalar extension $H$ of a Hopf algebra $T$ (with bijective antipode) by a braided commutative Yetter-Drinfeld module algebra $A$ we exhibit a Hopf $A$-algebroid with total algebra $H$ and with a balancing subalgebra given by a specified set of generators.
Clearly, every Hopf algebroid over a commutative base is both a Lu-Hopf algebroid and a Hopf algebroid with a balancing subalgebra, namely $\mathcal{B}=H\otimes_\genfd H$.

Lu~\cite{lu} exhibits an example which she calls a coarse Hopf algebroid, nowadays often called a {\bf minimal Hopf algebroid}. Given a unital associative algebra, the tensor product algebra $A\otimes A^{op}$ carries the structure of Hopf algebroid with source map $\alpha(a) = a\otimes 1$, target map $\beta(b) = 1\otimes b$, comultiplication $\Delta(a\otimes b) = (a\otimes 1)\otimes_A(1\otimes b)$, counit $\epsilon(a\otimes b) = a b$ and antipode $\tau(a\otimes b) = b\otimes a$. It has the balancing subalgebra $\mathcal{B} = (A\otimes\genfd)\otimes_A(\genfd\otimes A^{op})$.

\section{Scalar extension Hopf algebroids}
\label{sec:scalar}

\subsection{Scalar extensions, elements $R(a)$ and section $\gamma$}
Given any associative $\genfd$-algebra $A$ equipped with a left Hopf action $\blacktriangleright$ of a bialgebra $T$, vector space $A\otimes T$
carries a structure of a unital associative $\genfd$-algebra with 
multiplication bilinearly extending formula
$(a\otimes t)(a'\otimes t') = \sum a (t_{(1)}\blacktriangleright a')\otimes t_{(2)} t'$ and with unit $1_A\otimes 1_T$. This algebra is called the {\bf smash product algebra} (\cite{montg}) and denoted $A\sharp T$. It comes along with canonical algebra monomorphisms $A\cong A\otimes\genfd\hookrightarrow A\sharp T$ and $T\cong \genfd\otimes T\hookrightarrow A\sharp T$. The images of these two embeddings will be denoted $A\sharp 1$ and $1\sharp T$.

Let $T$ be a Hopf $\genfd$-algebra with a comultiplication 
$\Delta_T:T\to T\otimes_\genfd T$ and a bijective antipode $S$.
A {\bf braided-commutative left-right Yetter-Drinfeld 
$T$-module algebra} $A$ is a unital associative algebra
with a left $T$-action $\blacktriangleright:T\otimes A\to A$ 
which is Hopf in the sense that
$$
t\blacktriangleright (ab) = 
(t_{(1)}\blacktriangleright a)(t_{(2)}\blacktriangleright b),\,\,\,\,\,
t\blacktriangleright 1_A = \epsilon(t) 1_A,
$$
and a right $T$-coaction $a\mapsto a_{[0]}\otimes a_{[1]}$
which is morphism of algebras $A\to A\otimes T^{\mathrm{op}}$
~(see~\cite{BrzMilitaru}),
satisfying the left-right {\bf Yetter-Drinfeld condition}
\begin{equation}\label{eq:YD}
  (t_{(1)} \btr a_{[0]}) \otimes (t_{(2)} a_{[1]}) 
= (t_{(2)} \btr a)_{[0]} \otimes  (t_{(2)} \blacktriangleright a)_{[1]} t_{(1)},
\,\,\,\,\,\,\,\forall t \in T, \forall a\in A
\end{equation}
and the {\bf braided commutativity}
\begin{equation}\label{eq:bc}
  x_{[0]}(x_{[1]}\blacktriangleright a) = a x,\,\,\,\,\mbox{for }\,\,\,\,\,\,\forall a,x\in A.
\end{equation}
\begin{lemma}\label{lem:bcalt}
Braided commutativity condition is equivalent to the condition 
\begin{equation}\label{eq:bcalt}
  (S d_{[1]})\blacktriangleright a) d_{[0]} = d a,\,\,\,\,\forall d,a\in A.
\end{equation}
\end{lemma}
\begin{proof} This is rather standard. Assuming braided commutativity~(\ref{eq:bc}),
$$\begin{array}{lcl}
d a &=& d_{[0]} ((d_{[1]} S d_{[2]})\blacktriangleright a) \\
&=& d_{[0][0]}(d_{[0][1]}\blacktriangleright ((S d_{[1]})\blacktriangleright a))\\
&=& ((S d_{[1]})\blacktriangleright a) d_{[0]}.
\end{array}$$
  In other direction, assuming~(\ref{eq:bcalt})
  $$\begin{array}{lcl}
    x_{[0]}(x_{[1]}\blacktriangleright a)&=&
    (S x_{[0][1]}\blacktriangleright(x_{[1]}\blacktriangleright a))x_{[0][0]}
    \\
    &=& (S x_{[1][1]}\blacktriangleright(x_{[1][2]}\blacktriangleright a))x_{[0]}
    \\
    &=& (\epsilon(x_{[1]})1_T\blacktriangleright a)x_{[0]}
    \\
    &=& a x.
  \end{array}$$
\end{proof}

If $A$ is in fact a braided commutative Yetter-Drinfeld algebra over $T$ then the smash product $H = A\sharp T$ is a total algebra of a 
Hopf $A$-algebroid called a scalar extension Hopf algebroid.
For a Lu-Hopf algebroid this is proven in~\cite{BrzMilitaru}, modifying slightly an earlier construction of Lu~\cite{lu}, Section~5, where instead of Yetter-Drinfeld modules closely related modules over Drinfeld double $D(H)$ are considered. Both works entail a circular argument in the proof that the antipode of the algebroid is an antihomomorphism, checking the property on rather trivial case of binary products of generators of the form $a\sharp 1$ and $1\sharp t$ only. Antihomomorphism property for products of general elements is checked in~\cite{stojicPhD}, assuming that $S$ is bijective. 
Lu-Brzezi\'nski-Militaru construction has been adapted to the symmetric Hopf algebroids of B\"ohm~\cite{bohmnew,bohmHbk,BohmSzlach,stojicPhD}.

The $A$-bimodule structure of $A\sharp H$ is determined by the 
source and target maps
\begin{equation}
  \alpha(a) = a\sharp 1,\qquad \beta(a) = a_{[0]}\sharp a_{[1]},
\end{equation}
and the comonoid structure of $A\sharp H$ is given by 
\begin{equation}
  \Delta_{A\sharp T}(a\sharp t) = (a\sharp t_{(1)}) \otimes_A (1\sharp t_{(2)}),
  \qquad
  \epsilon_{A\sharp T}(a\sharp t) = a\epsilon_T(t). 
\end{equation}
Finally, the antipode $\tau$ for the Lu-Hopf algebroid is (cf.~\cite{BrzMilitaru,lu})  
\begin{equation}\label{eq:tause}
\tau(a\sharp t) = S(t) S^2(a_{[1]}) \cdot a_{[0]}.
\end{equation}
We often identify $A\sharp 1 = \mathrm{Im}\,\alpha$ with $A$
and $1\sharp T$ with $T$. 
By Definition~\ref{eq:IA} 
the right ideal $I_A\subset H\otimes_A H$ is generated by the set 
of all elements of the form
\begin{equation}\label{eq:Ia}
I(a) := \beta(a)\otimes 1 - 1\otimes\alpha(a) = 
a_{[0]}\sharp a_{[1]}\, \otimes 1 - 1 \otimes a,\,\,\,\,\,\,\,\,
a\in A.
\end{equation}

There is also another set of generators $R(a)$ of $I_A$, more convenient for 
our analysis below. 

\begin{proposition}\label{prop:RagenIA}
In the case of scalar extension $H = A\sharp T$, right ideal $I_A\in H\times_A H$ is generated by the set of all elements of the form 
\begin{equation}\label{eq:Ra}
R(a) := a \,\otimes 1 - S a_{[1]}\,\otimes a_{[0]},
\,\,\,\,\,\,\,\,\,a\in A. 
\end{equation}
\end{proposition}
\begin{proof} 
In the notation~(\ref{eq:Ia}),
$$
I(a) = (a_{[0]}\sharp 1 - S a_{[0][1]}\otimes a_{[0][0]})(a_{[1]}\otimes 1) 
= R(a_{[0]}) (a_{[1]}\otimes 1).
$$

Notice that $a_{[0]}\in A$. On the other hand, 
$$
R(a) = (a_{[0][0]} \sharp a_{[0][1]}\otimes 1 - 1\otimes a_{[0]}) (S a_{[1]}\otimes 1)
= I(a_{[0]})(S a_{[1]}\otimes 1). 
$$
Therefore, the right ideal generated by $\{I(a)\mid a\in A\}$ and the
right ideal generated by all $R(a)$ coincide. 
\end{proof}

Given a linear basis $\hx_1,\ldots,\hx_{\mathrm{dim}\,\mathfrak{g}}$
of a finite dimensional Lie algebra $\mathfrak{g}$, references~\cite{tajronkov,tajron} introduce elements $R_\mu$ ($\mu = 1,\ldots,\mathrm{dim}\,\mathfrak{g}$) in a related example of a formally completed version of a scalar extension Hopf $U(\gg)$-algebroid. In Subsection~\ref{ssec:comp}, we observe that $R_\mu = R(\hx_\mu)$.
\vskip .1in

\noindent {\bf Lu's section for scalar extension bialgebroids. }
For any scalar extension $H=A\sharp T$,
J-H.~Lu \cite{lu} exhibits
a section $\gamma : H\otimes_A H\to H\otimes_\genfd H$ 
by the unique $\genfd$-linear extension of the formula

\begin{equation}\label{eq:gammasmash}
\gamma : h \otimes_A (a\sharp t) \mapsto 
\beta(a) h \otimes_\genfd (1\sharp t),\,\,\,\,h\in H, a\in A, t\in T.
\end{equation}
Section $\gamma$ is well defined by~(\ref{eq:gammasmash}), 
namely on the generators
$$
\beta(b)h\otimes (c\sharp t) - h\otimes \alpha(b)(c\sharp t)
$$
of the ideal $I_A$
the formula evaluates to
$\beta(c)\beta(b)h\otimes(1\sharp t) -
\beta(b c)h\otimes(1\sharp t) = 0$.
Linear map $\gamma$ is a section of the projection 
$\pi : H\otimes_\genfd H\to H\otimes_A H$ because 
$h\otimes_A (a\sharp t) = h\otimes_A \alpha(a)(1\sharp t) = \beta(a) h\otimes_A (1\sharp t) = (\pi\circ\gamma)(h\otimes_A(a\sharp t))$. 

In particular, formula~(\ref{eq:gammasmash}) gives
\begin{equation}\label{eq:gammaDelta}
(\gamma\circ\Delta)(a\sharp t) = (a\sharp t_{(1)})\otimes_\genfd (1\sharp t_{(2)}).
\end{equation}

\subsection{Subalgebra $W\subset H\otimes H$ where $H=A\sharp T$ is a scalar extension Hopf $A$-algebroid}
\label{ssec:W}

\begin{notation} Let $T$ be a Hopf algebra and $A$ a braided commutative
algebra in the category of left-right Yetter-Drinfeld $T$-modules.
For a scalar extension $A\sharp T$ let
$W\subset (A\sharp T)\otimes (A\sharp T)$ 
be the smallest unital subalgebra such that all elements of 
the form $a\otimes 1$ and all elements of the form $S a_{[1]}\otimes a_{[0]}$
(where $a\in A\cong A\sharp 1\subset A\sharp T$) are in $W$. Let $W^+$
be the two sided ideal in $W$ generated by all elements of the
form $R(a) = a\otimes 1 - S a_{[1]}\otimes a_{[0]}$ 
where $a\in A$ (compare~(\ref{eq:Ra})).

Let $W_0^+\subset W$ be the linear subspace of $W$ spanned
by all elements of the form  $(x\otimes 1 - S x_{[1]}\otimes x_{[0]})(x'\otimes 1)$ where $x,x'\in A\cong A\sharp 1$. 
\end{notation}

We now formulate two lemmas which together imply $W_0^+ = W^+$.\vskip .1in 

\begin{lemma}\label{lem:R1}
For $x,z\in A$ we have $(x\otimes 1)(z\otimes 1 - S z_{[1]}\otimes z_{[0]})\in W_0^+$.
\end{lemma}

\begin{proof} Multiplying out, and using 
$x S(t) = S(t_{(1)}) t_{(2)} x S (t_{(3)}) 
= S (t_{(1)}) (t_{(2)}\blacktriangleright x)$ for $x\in A$, $t\in T$, we obtain
$$\begin{array}{lcl}
x z\otimes 1 - x S (z_{[1]})\otimes z_{[0]} &=&
x z \otimes 1 - S (z_{[1]}) (z_{[2]}\blacktriangleright x)\otimes z_{[0]}\\ 
&=& \mathrm{by\ braided\ commutativity} \\
&=& z_{[0]} (z_{[1]}\blacktriangleright x) \otimes 1 -
S (z_{[1]}) (z_{[2]}\blacktriangleright x)\otimes z_{[0]}\\
&=& (z_{[0]}\otimes 1 - S (z_{[0][1]}))\otimes z_{[0][0]}) 
(z_{[2]}\blacktriangleright x\otimes 1)
\end{array}$$
and the right hand side is clearly in $W_0^+$ as claimed.
\end{proof}

\begin{lemma} \label{lem:R2} 
$R(x)R(z) = (x\otimes 1 - S x_{[1]}\otimes x_{[0]})(z\otimes 1 - S z_{[1]}\otimes z_{[0]})\in W_0^+$.
\end{lemma}

\begin{proof} Since $x\mapsto x_{[1]}\otimes x_{[0]}$ 
is a morphism of algebras $A\to T^{\mathrm{op}}\otimes_\genfd A$
and $T\otimes_\genfd A\hookrightarrow A\sharp T\otimes_\genfd A\sharp T 
= H\otimes_\genfd H$
inclusion of algebras, we conclude that
$x\mapsto S x_{[1]}\otimes x_{[0]}$ is a morphism of algebras 
$A\to H\otimes_\genfd H$ 
(with respect to the componentwise multiplication in $H\otimes_\genfd H$). 
Therefore,

$$
(x\otimes 1 - S x_{[1]}\otimes x_{[0]})(z\otimes 1 
- S z_{[1]}\otimes z_{[0]}) = 
$$
$$
= 
(x\otimes 1 - S x_{[1]}\otimes x_{[0]})(z\otimes 1)
+ x S z_{[1]}\otimes z_{[0]} - 
x z\otimes 1 + x z\otimes 1 - S (x z)_{[1]}\otimes (x z)_{[0]}
$$
$$
= (x\otimes 1 - S x_{[1]}\otimes x_{[0]})(z\otimes 1) +
(-x\otimes 1) (z\otimes 1 
- S z_{[1]}\otimes z_{[0]}) + (x z\otimes 1 - S (x z)_{[1]}\otimes (x z)_{[0]}).
$$
The first and the third summands on the right hand side 
are manifestly in $W_0^+$ while for the second summand
we apply Lemma~\ref{lem:R1}.
\end{proof}

\begin{corollary} \begin{enumerate}[(i)] $\forall x,z\in A$, $(S x_{[1]}\otimes x_{[0]})(z\otimes 1-S z_{[1]}\otimes z_{[0]})\in W_0^+$,

\item $W^+_0$ is a two-sided ideal in $W$,

\item $W^+ = W_0^+$. 
\end{enumerate}
\end{corollary}

\begin{proof} (i) follows by subtracting the expression $(x\otimes 1 - S x_{[1]}\otimes x_{[0]})(z\otimes 1-S z_{[1]}\otimes z_{[0]})$ which is in $W^+_0$ by Lemma~\ref{lem:R2} from the expression 
  $(x\otimes 1)(z\otimes 1 - S z_{[1]}\otimes z_{[0]})$ which is in $W^+_0$ by Lemma~\ref{lem:R1}.

(ii) $W^+_0$ is a right ideal: by its definition, 
we can multiply by $z\otimes 1$ from the right;
this together with the assertion of Lemma~\ref{lem:R2} 
implies that we can also multiply
by $S z_{[1]}\otimes z_{[0]}$ from the right.

(ii) $W^+_0$ is a left ideal: using Lemma~\ref{lem:R1}, 
$(x\otimes 1) R(z) (x'\otimes 1) \in W_0^+ (x'\otimes 1)$
which is in $W_0^+$ because it is a right ideal. 
Combining with Lemma~\ref{lem:R2} we also conclude that 
$(S x_{[1]}\otimes x_{[0]})R(z)(x\otimes 1) \in W_0^+$.

For (iii) notice first that, trivially, $W_0^+\subset W^+$. 
For the converse inclusion, $W^+\subset W_0^+$, it is sufficient 
to observe that $R(a)\in W_0^+$, apply
(ii) and the definition of $W^+$.
\end{proof}

\begin{theorem} 
$\mu (id\otimes_\genfd\tau) W^+ = \{0\}$.
\end{theorem}

\begin{proof} 
By Corollary (iii) $W^+ = W_0^+$, which 
is the span of the elements of the form
$$(x\otimes 1 - S x_{[1]}\otimes x_{[0]})(z\otimes 1),\,\,\,\,\mathrm{where\ }\,
x,z\in A.$$
Taking the standard formula for the antipode for the scalar extensions~(\ref{eq:tause}),
$\tau(a\sharp t) = S(t) S^2(a_{[1]}) \cdot a_{[0]}$,
we can now compute $\mu (id\otimes\tau)$ on such an element as
$$
x z - S (x_{[2]}) z S^2(x_{[1]}) x_{[0]} 
= x z - ((S x_{[1]})\blacktriangleright z) x_{[0]}
= 0,
$$
by the braided commutativity. 
\end{proof}
\vskip .2in

\subsection{Subalgebra $\mathcal{B}$ and 
two-sided ideal $\mathcal{B}^+\subset\mathcal{B}$}
\label{ssec:B}

In this section, 
we want to show that every scalar extension Lu-Hopf algebroid 
$H = A\sharp T$ is also a Hopf algebroid with a 
(carefully chosen) balancing subalgebra $\mathcal{B}$.

Using the inclusion $T\otimes_\genfd T
\hookrightarrow A\sharp T\otimes_\genfd A\sharp T$,
we identify the image of the coproduct 
$\Delta_T:T\to T\otimes_\genfd T$ of the Hopf algebra $T$ 
with a subalgebra of $H\otimes_\genfd H$ which will be denoted 
by $\Delta_T(T)$. 

\begin{definition}\label{def:B}
Subalgebra $\mathcal{B}\subset A \sharp T \otimes_\genfd A \sharp T$ 
is the subalgebra generated by $W$ and $\Delta_T(T)$ 
or, equivalently, by the set
$$\left\{ a \otimes 1 ,\  Sa_{[1]}\otimes a_{[0]} \,\mid\,a\in A \right\} \cup \Delta_T(T).$$ 

The elements of this set are called the {\bf distinguished generators of} $\mathcal{B}$. Recall now elements $R(a)\in I_A\cap\mathcal{B}$
defined by formula~(\ref{eq:Ra}).
Let $\mathcal{B}^+$ be the two-sided ideal in $\mathcal{B}$ 
generated by the subset 
$$\left\{ R(a) \,\mid\,a\in A\right\} = \left\{ a\otimes 1   - \ Sa_{[1]}\otimes a_{[0]} \, \mid \,a\in A \right\}\subset\mathcal{B},$$
whose elements $R(a)$ are called the {\bf distinguished generators of} $\mathcal{B}^+$.
\end{definition}

\begin{theorem}\label{th:scalarbialgM} Suppose $H = A\sharp T$ is a
  smash product, where $T$ is a Hopf algebra with bijective antipode and $A$
  a braided-commutative Yetter-Drinfeld algebra over $A$; in other words,
  $A\sharp T$ is the underlying algebra of usual scalar extension bialgebroid.
  Let $\mathcal{B}$ and $\mathcal{B}^+$ be as in Definition~\ref{def:B}. 
\begin{enumerate}[(i)]
\item $\mathcal{B}^+ = I_A\cap\mathcal{B}$.

\item (C3MI) holds: $I_A\cap\mathcal{B}$ is a two-sided ideal
in $\mathcal{B}$. 

\item (C3Ma) holds: 
$\mathrm{Im}\,\Delta\subset\mathcal{B}/(I_A\cap\mathcal{B})$

\item The scalar extension $A\sharp T$ is a bialgebroid with the balancing subalgebra $\mathcal{B}$. 

\item $\mathcal{B}\subset \pi^{-1}(H\times_A H)$.

\item Inclusion from (v) induces an inclusion of algebras 
  $\mathcal{B}/\mathcal{B}_+\hookrightarrow H\times_A H$ whose image is
  $\Delta_L(A \sharp T) \subset A \sharp T \otimes_A A \sharp T$.
  On generators this isomorphism onto the image is given by $\genfd$-linear extension of the correspondence
  $a\otimes 1\mapsto a\sharp 1 \otimes_A 1$,
  $S a_{[1]}\otimes a_{[0]}\mapsto a\sharp 1\otimes_A 1$ for $a\in A$,
  and $\Delta_T(t) = t_{(1)}\otimes t_{(2)}\mapsto 1\sharp t_{(1)}\otimes_A 1\sharp t_{(2)}$ for $t\in T$.
\end{enumerate}
\end{theorem}

\begin{proof}
(i) follows immediately from Definition~\ref{def:B} 
  for $\mathcal{B}^+$ and Proposition~\ref{prop:RagenIA}.
  
  (ii) follows from (i) and the definition of $\mathcal{B}^+$.
  (iii) is an immediate check knowing the generators and (i).

For (iv) use (ii), (iii), Proposition~\ref{prop:equivifBexists}
and the fact that the (C3b) is known for scalar extensions~\cite{BrzMilitaru,lu}.

For (v) we have two proofs. One is to check for each distinguished generator separately that it belongs to $\pi^{-1}(H\times_A H)$. Another is to write $I(a)$ in terms of $R(a)$, use that $R(a)\in I_A\cap \mathcal{B}$ and hence $b R(a)\in b\mathcal{B}^+\subset \mathcal{B}^+$ because $R(a)\in\mathcal{B}^+$ and the latter is a two-sided ideal.

(vi) Clearly, $a\otimes 1 - S a_{[1]}\otimes a_{[0]}
= R(a) \in\mathcal{B}\cap I_A\mapsto I_A$, hence the values on $a\otimes 1$ and $S a_{[1]}\otimes a_{[0]}$ are the same. 
\end{proof}

\begin{lemma}\label{lem:key1}
   Let $\sum_{\rho_i} K_i^{\rho_i} \otimes L_i^{\rho_i} \in\left\{ x\otimes 1, \ Sx_{[1]}\otimes x_{[0]} \,|\,x\in A\right\}\cup \{f_{(1)}\otimes f_{(2)}\,|\,f\in T\}$ be a distinguished generator. Then for any $a\in A$,  
  $$\sum_{\rho_i} K_i^{\rho_i} \cdot (a \sharp 1) \cdot \tau(L_i^{\rho_i}) \in A \sharp 1.$$
\end{lemma}

\begin{proof}
  We inspect the claim case by case as follows.
\begin{itemize}
\item[(a)] If $\sum_{\rho_i} K_i^{\rho_i} \otimes L_i^{\rho_i} \in \left\{ x\otimes 1, \ Sx_{[1]}\otimes x_{[0]} \,|\,x\in A\right\}$, then by~(\ref{eq:tause})
\begin{equation}\label{eq:kZtL}
\sum_{\rho_i}K_i^{\rho_i} \cdot (a \sharp 1) \cdot \tau(L_i^{\rho_i}) = \sum_{\rho_i}K_i^{\rho_i} \cdot a\sharp 1 \cdot S^2(L_{i{[1]}}^{\rho_i}) \cdot L_{i[0]}^{\rho_i}.
\end{equation}

The dot product sign $\cdot$ denotes here the multiplication
in the smash product $A\sharp T$; $A$ is identified
there with $A\sharp 1$ and $T$ with $1\sharp T$.     
There are now two subcases, (a1) and (a2).
\begin{itemize}
\item[(a1)]
     For $\sum_{\rho_i} K_i^{\rho_i} \otimes L_i^{\rho_i} = x\otimes 1$, 
(\ref{eq:kZtL}) equals $xa = xa\sharp 1$ which is in $A\sharp 1$.
     
\item[(a2)]
 For $\sum_{\rho_i}K_i^{\rho_i} \otimes L_i^{\rho_i} = Sx_{[1]} \otimes x_{[0]}$, 
(\ref{eq:kZtL}) equals 
$Sx_{[2]} \cdot a \cdot S^2(x_{[1]}) \cdot x_{[0]} = ((Sx)_{[1]} \btr a) \cdot (Sx)_{[2]} S((Sx)_{[3]}) \cdot x_{[0]}  = ((Sx)_{[1]} \btr a) \cdot x_{[0]} = xa$ 
which is in $A \sharp 1$.
\end{itemize}     
\item[(b)] If $\sum_{\rho_i} K_i^{\rho_i} \otimes L_i^{\rho_i} = f_{(1)} \otimes f_{(2)}$, $f\in T$, then $\sum_{\rho_i} K_i^{\rho_i} \cdot a\sharp 1\cdot\tau(L_i^{\rho_i}) =  f_{(1)} \cdot a\sharp 1 \cdot S(f_{(2)}) =  (f \btr a)\sharp 1$ which is again in $A \sharp 1$.
 \end{itemize}
Therefore the claim $\sum_{\rho_i}K_i^{\rho_i}\cdot a\sharp 1 \cdot\tau(L_i^{\rho_i}) \in A \sharp 1$ follows.\end{proof}

\begin{lemma}\label{lem:key2}
   Let $U$ be a product of finitely many distinguished generators of $\mathcal{B}$. Then $$\mu(\id \otimes \tau)(U) \in  A \sharp 1.$$
\end{lemma}
\begin{proof} Let $U = (\sum_{\rho_1} K_1^{\rho_1} \otimes L_1^{\rho_1}) \cdots (\sum_{\rho_n} K_n^{\rho_n} \otimes L_n^{\rho_n})$. 
The antipode $\tau$ is an algebra antihomomorphism, hence
\begin{equation}\label{eq:Eexpanded}
  \mu(\id \otimes \tau)(U) = \sum_{\rho_1,\ldots,\rho_n}
  K_1^{\rho_1} K_2^{\rho_2} \cdots K_n^{\rho_n} \tau(L_n^{\rho_n}) \cdots \tau(L_1^{\rho_1}).
\end{equation}

We prove that
$$\sum_{\rho_{n-p},\ldots,\rho_n}K_{n-p}^{\rho_{n-p}}\cdots K_n^{\rho_n} \tau(L_n^{\rho_n})\cdots \tau(L_{n-p}^{\rho_{n-p}})\in A\sharp 1,$$
by induction on $p$ where $0\leq p\leq n-1$, the assertion of the lemma is then the case $p = n-1$. For the base of induction, $p=0$, this is the identity $\sum_{\rho_n} K_n^{\rho_n} \tau(L_n^{\rho_n}) \in A\sharp 1$ which follows from Lemma~\ref{lem:key1} when $a = 1$. The step of the induction on $p$ is clearly also a special case of Lemma~\ref{lem:key1}. 

By~(\ref{eq:Eexpanded}) it follows that $\mu(\id\otimes_\genfd \tau)(U) \in A \sharp 1$; in other words, $\mu(\id\otimes_\genfd \tau)(U)$ is of the form $d \sharp 1$ for some $d\in A$.
\end{proof}

\begin{theorem}\label{th:scalarM}
Let $\tau:A\sharp T\to A\sharp T$ given by the formula~(\ref{eq:tause})
be the antipode of the scalar extension as a Lu-Hopf algebroid.
Then
\begin{enumerate}[(i)]
\item $\mu( \id \otimes \tau) \mathcal{B}^+ = \left\{ 0 \right\} $. 
\item $\tau$ makes the corresponding $A$-bialgebroid with a balancing subalgebra
from Theorem~\ref{th:scalarM}
into a Hopf $A$-algebroid with a balancing subalgebra.
\end{enumerate}
\end{theorem}

\begin{proof} (i)
A general element of $\mathcal{B}^+$ is a linear combination of the
elements of the form
\begin{equation}\label{eq:genMNxKL}
  \prod_{j=1}^m \sum_{\sigma_j} M_{j}^{\sigma_j} \otimes N_{j}^{\sigma_j} \cdot \left( x \otimes 1 - Sx_{[1]} \otimes x_{[0]} \right) \cdot \prod_{k=1}^{n} \sum_{\rho_k} K_{k}^{\rho_k} \otimes L_{k}^{\rho_k},
\end{equation}
where $\sum_{\sigma_j} M_{j}^{\sigma_j} \otimes N_{j}^{\sigma_j}$, $\sum_{\rho_k}K_{k}^{\rho_k} \otimes L_{k}^{\rho_k}$
are some distinguished generators of $\mathcal{B}$, and the midd\-le factor
$x\otimes 1 - S x_{[1]}\otimes x_{[0]}$ is some distinguished generator in $\mathcal{B}^+$. 
Notice that 
$M_{j}^{\sigma_j}$, $N_{j}^{\sigma_j}$, $K_{k}^{\rho_k}$, $L_{k}^{\rho_k} \in A \sharp 1 \cup 1 \sharp T$ and $x \in A$.

By the linearity of $\mu( \id \otimes \tau)$ it is sufficient to prove
the assertion for one element of the form above. 
Rewrite~(\ref{eq:genMNxKL}) as,
\begin{equation}\label{eq:gen2}
  \prod_{j,k}\sum_{\sigma_j,\rho_k}(M_j^{\sigma_j} \cdot x\sharp 1  \cdot K_k^{\rho_k} )\otimes (N_j^{\sigma_j} \cdot L_k^{\rho_k}) - (M_j^{\sigma_j} \cdot 1 \sharp Sx_{[1]} \cdot K_k^{\rho_k} ) \otimes (N_j^{\sigma_j} \cdot x_{[0]} \sharp 1 \cdot L_k^{\rho_k} ).
  \end{equation}
By Lemma~\ref{lem:key2}, we can define $d\in A$ by  
\begin{equation}\label{eq:e}
d \sharp 1 :=  K_1^{\rho_1}\cdots K_n^{\rho_n} \cdot \tau(L_n^{\rho_n})\cdots\tau(L_1^{\rho_1}) \in A\sharp 1.
\end{equation}
We apply map $\mu( \id \otimes \tau)$ to~(\ref{eq:gen2}) and substitute~(\ref{eq:e}). Notice that $\tau$, being an antihomomorphism, reverses the order. Thus we need the vanishing of
\begin{equation}\label{eq:MxetN}
\sum_{\sigma_1\ldots\sigma_m} M_1^{\sigma_1}\cdots M_m^{\sigma_m}
\left(x\sharp 1  \cdot d\sharp 1
-  1 \sharp Sx_{[1]} \cdot d\sharp 1 \cdot  \tau(x_{[0]} \sharp 1)\right)
\tau(N_m^{\sigma_m})\cdots\tau(N_1^{\sigma_1}).
\end{equation}
      
     Therefore to finish the proof of the assertion (i) it is 
sufficient to show that for all $x$, $d \in A$ we have     
     $$ x\sharp 1  \cdot d \sharp 1 =  1 \sharp Sx_{[1]} \cdot d \sharp 1 \cdot  \tau( x_{[0]} \sharp 1 ),$$    
   where, by the formula for the antipode~(\ref{eq:tause}),  $\tau(x_{[0]} \sharp 1) = S^2(x_{[1]}) \cdot x_{[0]}$.    
     
This amounts to showing 
$$xd \sharp 1 = \left((Sx)_{[1]} \btr d \right) \sharp \,  (Sx)_{[2]} \,  S((Sx)_{[3]}) \cdot x_{[0]} \sharp 1 ,$$
that is,
$$xd \sharp 1 = \left((Sx)_{[1]} \btr d\right) x_{[0]} \sharp 1,$$
which is by Lemma~\ref{lem:bcalt} an expression of 
the braided commutativity~(\ref{eq:bc}) for $A$. 
Therefore, (i) is proven.

For the part (ii), according to Theorem~\ref{th:scalarbialgM}, part (iv), 
it remains only to check the axioms for the antipode.
The antipode requirements~(\ref{eq:mtaubetaalpha}) and~(\ref{eq:ml}) 
have the same content 
as in the case of Lu-Hopf algebroid definition hence they are true.
Now, thanks to (ii) the left-hand side of the equation~(\ref{eq:mr}),
that is, $\mu( \id \otimes_A \tau) \Delta$, 
does not depend on the representatives of $\Delta(a t)
= (a\sharp t_{(1)})\otimes (1\sharp t_{(2)})$ in $H\otimes_\genfd H$
where $a\in A$ and $t\in T$. So we need to show that  
$$
(a \sharp t_{(1)}) (S^2(t_{(2)})_{[1]}\cdot (t_{(2)})_{[0]}) = a\sharp t,
$$ 
which boils down to the same computation as
for the Lu's choice of $\gamma$, see~(\ref{eq:gammaDelta}).
Our result is stronger only in the sense that we allow for an additional freedom in $\mathcal{B}_+$ and that 
$\mathcal{B}$ is a balancing subalgebra in the bialgebroid sense.
\end{proof}

\subsection{Comparison with the examples of Meljanac}
\label{ssec:comp}

S.\ Meljanac has devised his method~\cite{tajronkov,tajron} 
to the study topological Hopf algebroids related to
a Lie algebra $\gg_\kappa$ with 
the universal enveloping algebra $U(\gg_\kappa)$ 
in physics literature called the $\kappa$-Minkowski space.
Some extensions of this Hopf algebroid
(including some symmetries into the algebra)
from the point of view of Lu-Hopf algebroid
have been studied in~\cite{kappaPLB} in
an informal style of mathematical physics
and, in just slightly more mathematical treatment, in~\cite{JLZSMWa}. 
Works~\cite{tajronkov,tajron} made it clear that their
construction applies to any finite dimensional Lie algebra
$\gg$ in characteristic zero. We comment below on 
how our construction of $\mathcal{B}$ 
relates to theirs for general $\gg$.
As stated in the introduction, we neglect here the issues
related to the adaptation of the notion of Hopf algebroid to 
the completed tensor products~\cite{halgoid,stojicPhD}. 

We use the notation from~\cite{halgoid}.
Generators of the Lie algebra $\gg$ 
are denoted $\hat{x}_1,\ldots,\hat{x}_n$ 
with commutators $[\hat{x}_\mu,\hat{x}_\nu] = C_{\mu\nu}^\lambda\hat{x}_\lambda$
and the generators of
the symmetric algebra of the dual $S(\gg^*)$ 
by $\partial^1,\ldots,\partial^n$. The completed dual 
$T = \hat{S}(\gg^*)$ is a topological Hopf algebra, namely the coproduct $\Delta_T:\hat{S}(\gg^*)\to\hat{S}(\gg^*)\hat\otimes\hat{S}(\gg^*)$
may be identified with the dual (transpose) map
to the multiplication $U(\gg)\to U(\gg)\otimes U(\gg)$.
The identification is made with help of 
the symmetrization map $S(\gg)\cong U(\gg)$,
which is an isomorphism of coalgebras~\cite{ldWeyl,halgoid}
and its dual isomorphism of algebras $\hat{S}(\gg^*)\cong U(\gg)^*$.
Now $A = U(\gg)$ becomes a braided-commutative
Yetter-Drinfeld module algebra over $T$ (internally in a symmetric category of filtered-cofiltered vector spaces,~\cite{stojicPhD}). Regarding that $T$ is a formal dual of $U(\gg)$, the Heisenberg double $U(\gg)\sharp T$~\cite{halgoid,heisd,stojicPhD} can be either produced as a usual smash product where $T$ is equipped with a right Hopf action of a Hopf algebra $U(\gg)$ (as in~\cite{halgoid}, where however, for a bialgebroid, an additional completion on the smash product has been performed at a later stage) or a smash product in which $T$ is understood as a topological Hopf algebra and $U(\gg)$ is equipped with an internal left Hopf action of $T$. The latter interpretation supplies the internal version of a scalar extension Hopf algebroid~\cite{stojicPhD}.

If $R(a)$ is defined by~(\ref{eq:Ra}), then
$\hat{x}_{\mu[0]}\otimes\hat{x}_{\mu[1]} = \hat{x}_\sigma\otimes(\mathcal{O}^{-1})^\sigma_\mu$ where $\mathcal{O}^\sigma_\tau,(\mathcal{O}^{-1})^\sigma_\tau$ are certain elements in $T$ (see~\cite{halgoid} for the definition and properties),
$\mathcal{O}^{-1}$ is a matrix inverse of $\mathcal{O}$,
$\Delta_T\mathcal{O}^\mu_\nu = \mathcal{O}^\sigma_\nu\otimes\mathcal{O}^\mu_\sigma$ and $S(\mathcal{O}^{-1})^\mu_\nu = \mathcal{O}^\mu_\nu$.
Thus we obtain,
\begin{equation}\label{eq:Rxmu}
  R(\hat{x}_\mu)= \hat{x}_\mu\otimes 1 - S\hat{x}_{\mu[1]}\otimes\hat{x}_{\mu[0]}
   = \hat{x}_\mu \otimes 1 - \mathcal{O}^\sigma_\mu\otimes\hat{x}_\sigma.
\end{equation}
We observe that $R(\hat{x}_\mu)$ is identical to $R_\mu$ of~\cite{tajron,tajronkov}.
Using identities $[\mathcal{O}^\sigma_\mu,\hat{x}_\nu] = C^\rho_{\mu\nu} O^\sigma_\rho$ (formula~(17) in~\cite{halgoid})
and $C_{\mu\nu}^\tau\mathcal{O}^\lambda_\tau = C^\lambda_{\rho\sigma}\mathcal{O}^\rho_\mu\mathcal{O}^\sigma_\nu$ (formula~(20) in~\cite{halgoid}), we obtain
\begin{equation}\label{eq:RmuRnu}
 [R(\hat{x}_\mu), R(\hat{x}_\nu)]  = C^\sigma_{\mu\nu} R(\hat{x}_\sigma),
\end{equation}
\begin{equation}\label{eq:xmuRnu}
    [\hat{x}_\mu\otimes 1,R(\hat{x}_\nu)] = C_{\mu\nu}^\lambda R(\hat{x}_\lambda)
\end{equation}
(generalizing Eq.~(32),(33) in arXiv version of~\cite{tajronkov}, (3.2),(3.3) in journal v.). Moreover, (\ref{eq:RmuRnu}) is the only relation among $R_\mu$-s, hence the subalgebra
in $H\otimes_\genfd H$ generated by $\{R_\mu\mid \mu =1, \ldots, n\}$,
is isomorphic to the universal enveloping algebra $U(\gg)$,
but with generators $R_\mu$ in place of $\hat{x}_\mu$.
Following~\cite{tajronkov}, denote this subalgebra by $U(R)$. 
The relation~\cite{tajronkov}~(\ref{eq:xmuRnu})
shows that the products of the form $r(u\sharp 1\otimes 1)$
where $r\in U(R)$ and   
$u\sharp 1\otimes 1 \in U(\gg)\sharp 1 \otimes_\genfd \genfd\subset
U(\gg)\sharp T\otimes_\genfd U(\gg)\sharp T$ 
span a subalgebra in $H\otimes_\genfd H$. 
This is precisely our subalgebra $W$ in this case.
However, the relations~(3.3),(3.4) in~\cite{tajronkov}
for $\gg_\kappa$ case,
and the generalizations~(\ref{eq:RmuRnu}), (\ref{eq:xmuRnu})
for general $\gg$,
used to show that $W$ is a subalgebra of $H\otimes_\genfd H$,
do not have a simple analogue for general scalar extension $A\sharp T$
(not of enveloping algebra type). 
It is also not clear what is the precise 
structure of the subalgebra generated by $R(a)$-s for all $a\in A$,
in general. 
On the other hand, our Hopf algebraic definition~(\ref{eq:Ra}) of
$R(a)$ and the corresponding definition of $W$ 
in the subsection~\ref{ssec:W} along with the lemmas
therein guarantee that such general $W$ is a subalgebra in 
$H\otimes_\genfd H$ in full generality.

The issues are more complicated when we pass from 
$W$ to $\mathcal{B}$.
In the enveloping algebra case, the subalgebra $\mathcal{B}$ 
(denoted $\hat{\mathcal{B}}$ in~\cite{tajronkov})
is defined in~\cite{tajronkov} rather simply as the subalgebra of
all elements of the form $\sum_i w_i \Delta_{\hat{S}(\gg^*)}(t_i)$
where $w_i\in W$, $t_i\in T = \hat{S}(\gg^*)$ are arbitrary
(the sums may be infinite, in an appropriate completion).
Equations~(3.3),(3.4) in~\cite{tajronkov}, can be abstracted and
generalized to an 
arbitrary finite dimensional Lie algebra as the following
proposition.
\begin{proposition} For general $\gg$, 
  \begin{equation}\label{eq:gen33}\begin{array}{l}
  \lbrack\Delta \partial^\mu,R(\hat{x}_\nu)\rbrack = 0 \\ 
  \lbrack\Delta\partial^\mu,\hat{x}_\nu\sharp 1 \otimes 1\rbrack
  \in\Delta_T(T).
\end{array}\end{equation}
\end{proposition}
Regarding that $\partial^\mu$ generate a dense subalgebra of $T$,
this implies immediately that
$\{\sum_i w_i\Delta_{\hat{S}(\gg^*)}(t)\mid w_i\in W, t_i\in T\}$
is a subalgebra of $H\otimes_\genfd H$, and that 
$\mathcal{B}$ has a very simple structure of 
all sums of products of the form: 
an element in $U(R)$ times 
an element in $A\sharp 1\otimes_\genfd \genfd
\subset H\otimes_\genfd H$
times an element of the form $\Delta_T(t_i)$ with $t_i \in T$.
We have exhibited above a similar structure
-- as a sum of products of elements from three subalgebras in 
this fixed order --
for general scalar extension Hopf algebroids. In this generality, 
$P_i$ do not commute with elements in $W$ and multiple 
products (e.g.\ of the form $w t w' t' w''$) of elements in $W$ and 
elements in $\Delta_T(T)$ may 
appear, as analysed in the subsection~\ref{ssec:B}.
Regarding that $\mu(\id\otimes_\genfd\tau)$ is not an antihomomorphism,
the multiple products bring the main difficulty in
our proof that the antipode $\tau$ is well defined 
(see Theorem~\ref{th:scalarM} (i)).

Analogous comparisons may be made for the ideal
$\mathcal{B}_+$ 
which is in~\cite{tajronkov} not defined as the intersection $I_A\cap\mathcal{B}$, but an equivalent description is given, constructing it 
in analogy to $\mathcal{B}$, but with the enveloping algebra $U(R)$
replaced by its ideal $U_+(R)\subset U(R)$ of elements which
are not degree $0$ in the standard filtration of the 
universal enveloping algebra. This explains the notation $\mathcal{B}_+$.
The commutation relations (3.1)-(3.4) in~\cite{tajronkov} imply
that such $\mathcal{B}_+$ is indeed a two-sided ideal in $\mathcal{B}$.

Our approach also differs from~\cite{tajronkov} in insisting
that the coproduct is still defined as taking values in $H\otimes_A H$
(rather than in $\mathcal{B}/(I_A\cap\mathcal{B})$ as an abstract algebra);
the two-sided ideal trick is used only to make sense of the requirement
and to check that the induced
map into $\mathcal{B}/(I_A\cap\mathcal{B})$ is a morphism of algebras. 
Moreover, they view $\mathcal{B}$ as an abstract algebra
constructed from its pieces $U(R)$, $A\sharp 1\otimes_\genfd\genfd$
and $\Delta_T(T)\cong T$.
In our approach, the coherently associative 
tensor product of bimodules $\otimes_A$ is used to formulate
the coassociativity of the coproduct as in the standard
definition of a bialgebroid~\cite{bohmHbk,BohmSzlach,BrzMilitaru,lu,xu}.
In~\cite{tajronkov} the coproduct is taking values in
$\mathcal{B}/(I_A\cap\mathcal{B})$ by definition
and, in the spirit of their viewpoint, the higher iterations
of the coproduct in subalgebras 
$$\mathcal{B}^{(j)}\subset 
H\otimes_\genfd H\otimes_\genfd \cdots\otimes_\genfd H\,\,\,\,
(j\,\,\,\mathrm{tensor}\,\,\,\mathrm{factors}),
$$
which define higher analogues of 
the subalgebra $\mathcal{B}^{(2)}:=\mathcal{B}\subset H\otimes_\genfd H$.
One also considers
the higher analogues $\mathcal{B}^{(j)}_+ = I_A^{(j)}\cap\mathcal{B}^{(j)}$
of $I_A\cap\mathcal{B}$ 
in order to deal with the (co)associativity issues.
For example, $\mathcal{B}^{(3)}$ is generated by all ordered products of the form $r\cdot(a\otimes 1\otimes 1) \cdot (\Delta\otimes\id)(\Delta(t))$
 where $a\in A$, $t\in T$ and $r$ belongs to the subalgebra generated by
$\{R(a)\otimes 1\mid a\in A\}\cup\{1\otimes R(b)\mid b\in A\}$.
 The right ideal $I_A^{(j)}$ is the smallest right ideal in $H^{\otimes k}$ containing right ideals $I_A\otimes H^{\otimes (j-2)}, H\otimes I_A \otimes H^{\otimes (j-3)},\ldots,H^{\otimes (j-2)}\otimes I_A$.
These are interesting structures, but in our view more cumbersome than
the familiar usage of the bimodule tensor product $\otimes_A$.

\section{Weak Hopf algebras}\label{sec:wha}

It is well-known that from the data of any weak Hopf algebra
one can construct a corresponding Lu-Hopf algebroid. 
Upon looking at our axiomatics, 
G.~B\"ohm has observed and sketched to us how to
construct a Hopf algebroid with balancing subalgebra 
from a weak Hopf algebra. 
We present her results in this section, starting with
a short review of weak Hopf algebras.

\subsection{Weak bialgebras, standard definitions}
\def\HH{\mathbb H}

A {\bf weak $\genfd$-bialgebra} $\HH$ (see~\cite{bohmnillszlachwHa})
is a tuple $(\HH,\mu,\eta,\Delta,\epsilon)$ where
$(A,\mu,\eta)$ is an associative unital $\genfd$-algebra, 
$(\HH,\Delta,\epsilon)$ is a coassociative counital $\genfd$-coalgebra,  
and the following compatibilities hold:
\begin{enumerate}[(i)]
\item $\Delta$ is multiplicative, $\Delta(a b) =
\Delta(a)\Delta(b)$ for all $a,b\in \HH$.
\item Weak multiplicativity of the counit: for all $x,y,z\in\HH$,
\begin{equation}\label{eq:wmult}
\epsilon(x y z) = \epsilon(x y_{(1)})\epsilon(y_{(2)} z),
\end{equation}
\begin{equation}\label{eq:wmultop}
\epsilon(x y z) = \epsilon(x y_{(2)})\epsilon(y_{(1)} z).
\end{equation}
If we assume (i), then it is elementary (see formulas~(4) and~(1) in~\cite{bohmCaenJans2011}) that~(\ref{eq:wmult}) and~(\ref{eq:wmultop})) are respectively equivalent to the conditions
\begin{equation}\label{eq:wmout}
  g\epsilon(1_{(2)}h)1_{(1)}=\epsilon(g_{(2)}h)g_{(1)}, \,\,\,\,\,\,\,\,\,\forall g,h\in\HH,
\end{equation}
\begin{equation}\label{eq:wmopout}
  g\epsilon(1_{(1)}h)1_{(2)}=\epsilon(g_{(1)}h)g_{(2)}, \,\,\,\,\,\,\,\,\,\forall g,h\in\HH.
\end{equation}

\item Weak comultiplicativity of the unit:

\begin{equation}\label{eq:Delta2}\begin{array}{l}
\Delta^{(2)} (1) = (\Delta(1) \otimes 1)(1\otimes \Delta(1))
\\
\Delta^{(2)} (1) = (1 \otimes\Delta(1))(\Delta(1) \otimes 1)
\end{array}\end{equation}
where we denoted $\Delta^{(2)} := (\id\otimes\Delta)\Delta = 
(\Delta\otimes\id)\Delta$. In Sweedler notation, 
\begin{equation}\label{eq:Delta2bis}
  1_{((1)}\otimes 1_{(2)}\otimes 1_{(3)} =
  1_{(1)}\otimes 1_{(2)}1_{(1)'}\otimes 1_{(2)'}=
  1_{(1)'}\otimes 1_{(1)}1_{(2)'}\otimes 1_{(2)}.
\end{equation}
\end{enumerate}

For every weak $\genfd$-bialgebra there are $\genfd$-linear maps 
$\Pi^L,\Pi^R:\HH\to\HH$
with properties $\Pi^R\Pi^R = \Pi^R$ and $\Pi^L\Pi^L = \Pi^L$
and defined by
$$
\Pi^L(x) := \epsilon(1_{(1)} x) 1_{(2)},\,\,\,\,
\Pi^R(x) := 1_{(1)}\epsilon(x 1_{(2)}).
$$
These expressions are met below in two of the axioms for the antipode
of a weak Hopf algebra. Less frequently, one also encounters idempotents
$\bar\Pi^L,\bar\Pi^R$ given by $$
\bar\Pi^L(x) := \epsilon(1_{(2)} x) 1_{(1)},\,\,\,\,
\bar\Pi^R(x) := 1_{(2)}\epsilon(x 1_{(1)}).
$$
Now
$$\begin{array}{lcl}
\epsilon(x z) = \epsilon(x 1 z) 
\stackrel{(\ref{eq:wmultop})}= \epsilon(x 1_{(2)})\epsilon(1_{(1)}z) 
&=& \epsilon(\epsilon(x 1_{(2)})1_{(1)}z)
= \epsilon(\Pi^R(x)z),
\\
&=& \epsilon(x \epsilon(1_{(1)}z)1_{(2)}) 
= \epsilon(x\Pi^L(z)).
\end{array}$$
The images of the idempotents $\Pi^R$ and $\Pi^L$, 
$$
\HH^R := \Pi^R(\HH),\,\,\,\,\,\,\,\,\HH^L = \Pi^L(\HH),
$$
are mutually dual as $\genfd$-linear spaces via
the canonical nondegenerate pairing $\HH^L\otimes\HH^R\to\genfd$ 
given by $(x,y) \mapsto \epsilon(y x)$.

The identities $\Pi^L(x\Pi^L(y)) = \Pi^L(x y)$ 
 and $\Pi^R(\Pi^R(x)y) = \Pi^R(x y)$ hold. 
Dually also $\Delta(\HH^L)\subset \HH\otimes\HH^L$, 
$\Delta(\HH^R)\subset \HH^R\otimes\HH$,
 and $\Delta(1)\in \HH^R\otimes \HH^L$.
 
\subsection{B\"ohm's recipes}

It is known that a weak Hopf algebra $\HH$ 
can be regarded as a Lu-Hopf algebroid over $A:=\Pi^L(\HH)$
where the source map $\alpha$ being the inclusion $\Pi^L(\HH)\subset\HH$ 
and where the target map is given by
\begin{equation}\label{eq:betawha}
\beta(a) = \bar\Pi^L(a) = \epsilon(1_{(2)}a) 1_{(1)}\,\,\,\,\,\mbox{for}\,\,\,\,a\in A,
\end{equation}
and the comultiplication $\Delta' = \pi\circ\Delta$ of the bialgebroid 
is the comultiplication $\Delta:\HH\to\HH\otimes_\genfd\HH$
of the weak Hopf algebra 
followed by the canonical
projection $\pi:\HH\otimes_\genfd\HH\to\HH\otimes_{\Pi^L(\HH)}\HH$.
\begin{lemma}\label{lem:difunit}
$1\otimes 1-\Delta(1)\in I_A$.
\end{lemma}
\begin{proof}
By definition, $I_A$ is generated by all expressions of the form
$$
I(h) := \beta(\Pi^L(h))\otimes 1 - 1\otimes\alpha(\Pi^L(h)),\,\,\,h\in\HH.
$$
Now $\Pi^L(h) = \epsilon(1_{(1)} h) 1_{(2)}$
hence by~(\ref{eq:betawha})
\begin{equation}\label{eq:betaPiL}
\begin{array}{lcl}\beta(\Pi^L(h)) &=&
\epsilon(1_{(2)}\epsilon(1_{(1')}h)1_{(2')})1_{(1)}\\ &=&
\epsilon(1_{(2)}1_{(2')})\epsilon(1_{(1')}h)1_{(1)}\\ &\stackrel{(\ref{eq:wmultop})}{=}&
\epsilon(1_{(2)}h)1_{(1)},\end{array}\end{equation}
where $1_{(1')}\otimes 1_{(2')}$ denotes another copy of $\Delta(1)$.
\begin{equation}\label{eq:Ihwha}\begin{array}{lcl}
I(h) &=& \epsilon(1_{(2)}h)1_{(1)} \otimes 1
- 1 \otimes \epsilon(1_{(1)} h) 1_{(2)}\\
&=&\bar\Pi^L(h)\otimes 1 - 1 \otimes\Pi^L(h).
\end{array}
\end{equation}

It is sufficient to prove that 
$1\otimes 1 -\Delta(1) = I(1_{(2)})(1_{(1)}\otimes 1)$ because the right hand side manifestly belongs to $I_A$. 
From~(\ref{eq:Ihwha}) we calculate
$$\begin{array}{lcl}
I(1_{(2')})(1_{(1')}\otimes 1) &=& \epsilon(1_{(2)}1_{(2')}) 1_{(1)}1_{(1')}
\otimes 1 - 1_{(1')} \otimes \epsilon(1_{(1)}1_{(2')})1_{(2)}
\\
&=& \epsilon((1\cdot 1)_{(2)}) (1\cdot 1)_{(1)} \otimes 1 
- 1_{(1)}\otimes \epsilon(1_{(2)}) 1_{(3)} 
\\
&=& 1\otimes 1 - \Delta(1),
\end{array}$$
where in the middle line the axioms~(\ref{eq:Delta2}) on $\Delta^{(2)}$ were used for the second summand.
\end{proof}

\begin{lemma}\label{lem:Igen}
  $\Delta(1)I(h) = 0$. 
\end{lemma}
\begin{proof} By~(\ref{eq:Ihwha}), $\Delta(1)I(h)=1_{(1)}\epsilon(1_{(2')}h)1_{(1')}\otimes 1_{(2)}-1_{(1)}\otimes 1_{(2)}\epsilon(1_{(1')}h)1_{(2')}$. 
$$\begin{array}{lcl}
1_{(1)}\epsilon(1_{(2')}h)1_{(1')}\otimes 1_{(2)}
&\stackrel{(\ref{eq:wmout})}=&
\epsilon(1_{(1)(2)}h)1_{(1)(1)}\otimes 1_{(2)}
\\ &=& 1_{(1)}\otimes \epsilon(1_{(2)(1)}h) 1_{(2)(2)}
\\
&\stackrel{(\ref{eq:wmopout})}=& 1_{(1)}\otimes 1_{(2)}\epsilon(1_{(1')}h)1_{(2')}
  \end{array}$$
Here we used (\ref{eq:wmout}) with $g=1_{(1)}$ and (\ref{eq:wmopout}) with $g=1_{(2)}$.
\end{proof}
\begin{corollary}
The right ideal $I_A$ coincides with the principal right ideal generated by $1\otimes 1 - \Delta(1)$.
\end{corollary}
\begin{proof}
  By Lemma~\ref{lem:difunit} element $1\otimes 1-\Delta(1)\in I_A$ and by Lemma~\ref{lem:Igen}, for every $h\in\HH$, $I(h) = (1\otimes 1 - \Delta(1))I(h)$, hence also $I_A\subset (1\otimes 1 -\Delta(1))\HH$. 
\end{proof}
\begin{theorem}
For a weak bialgebra $(\HH,\mu,\eta,\Delta,\epsilon)$, define the subalgebra 
$$\mathcal{B}:= \Delta(1)(\HH\otimes\HH)\Delta(1)\subset \HH\otimes\HH.$$ 
Then $(\HH,\mu,\Pi^L(\HH),\bar\Pi^L,\pi|_{\mathcal{B}}\circ\Delta,\Pi^L)$ is a left $\Pi^L(\HH)$-bialgebroid with balancing subalgebra $\mathcal{B}$.
\end{theorem}
\begin{proof}
It is clear that $\mathcal{B}$ is a subalgebra with unit $\Delta(1)$ 
and that $\operatorname{Im}\,\Delta\subset\mathcal{B}$.
By Lemma~\ref{lem:Igen}, the intersection $\mathcal{B}\cap I_A \subset \Delta(1)I_A = 0$ is the zero ideal of $\mathcal{B}$ hence (C3MI) holds and $\Delta':\HH\longrightarrow\HH\otimes_{\Pi^L(\HH)}\HH$
factorizes, indeed, through an algebra homomorphism 
$\HH\longrightarrow\mathcal{B}/(\mathcal{B}\cap I_A)$.
It is clear that $\Delta(h) = \Delta(1)\Delta(h)\Delta(1)\in\mathcal{B}$ for every $h\in\HH$, hence (C3Ma) holds. 
Then $\pi|_{\mathcal{B}}:\mathcal{B}\cong \mathcal{B}/(\mathcal{B}\cap I_A)$, hence as $\Delta$ is homomorphism, its corestriction $\Delta|^{\mathcal{B}}$ to $\mathcal{B}$ is homomorphism and, finally, the corestriction followed by the restriction of the projection $\pi|_{\mathcal{B}}\circ\Delta|^{\mathcal{B}}$ is a homomorphism. Other properties (e.g. that $(\HH,\Delta,\Pi^L)$ is a $\Pi^L(\HH)$-coring) are well known as they coincide with the axioms
of a left associative $A$-bialgebroid.
\end{proof}

\subsection{Antipode}

A weak $\genfd$-bialgebra $\HH$ is a {\bf weak Hopf algebra} 
if there is a $\genfd$-linear map $S:\HH\to\HH$ 
(which is then called an antipode)
such that for all $x\in\HH$
\begin{equation}\label{eq:wHAS1}
x_{(1)} S(x_{(2)}) = \epsilon(1_{(1)} x)1_{(2)},
\end{equation}
\begin{equation}\label{eq:wHAS2}
S(x_{(1)})x_{(2)} = 1_{(1)} \epsilon(x 1_{(2)}),
\end{equation}
\begin{equation}\label{eq:wHAS3}
S(x_{(1)})x_{(2)} S(x_{(3)}) = S(x)
\end{equation}
Notice that the right hand side of~(\ref{eq:wHAS1}) equals $\Pi^L(x)$
and the right hand side of~(\ref{eq:wHAS2}) equals $\Pi^R(x)$.
Suppose the antipode $S$ is bijective.
Set the antipode of the corresponding Hopf algebroid 
with a balancing subalgebra to be $\tau = S$. 
Since $I_A\cap\mathcal{B} = \{0\}$  
any $\genfd$-linear map vanishes on it; 
hence so does the map $\mu\circ(\id\otimes_\genfd\tau)$ of~(\ref{eq:tauvanish}). 
Axiom~(\ref{eq:mtaubetaalpha}) can be restated as
$(S\circ\beta\circ\Pi^L)(h) = \Pi^L(h)$. To show this identity,
notice that $(S\circ\beta\circ\Pi^L)(h) = S(\epsilon(1_{(2)}h)1_{(1)})$
by~(\ref{eq:betaPiL}) and then 
it is enough to quote $S(1_{(1)})\epsilon(1_{(2)}x) = \Pi^L(x)$,
which is the identity (2.24a) in~\cite{bohmnillszlachwHa}. Axiom~(\ref{eq:mr})
reads $h_{(1)}S(h_{(2)}) = \Pi^L(h)$ which is manifestly~(\ref{eq:wHAS1}).
Axiom~(\ref{eq:ml}) follows by calculation $S(h_{(1)})h_{(2)} = \beta(\Pi^L(S(h))) \stackrel{(\ref{eq:betaPiL})}= 1_{(1)}\epsilon(1_{(2)}S(h)) = 1_{(1)}\epsilon(h 1_{(2)})$, where the last equality is (2.23b) in~\cite{bohmnillszlachwHa}, proven using axiom~(\ref{eq:wHAS3}).

\section{Twisting by invertible 2-cocycles}

Ping Xu~\cite{xu} has generalized Drinfeld's procedure of twisting of bialgebras by invertible counital 2-cocycles to associative bialgebroids. Basic treatment involves several subtle points~\cite{xu} not appearing in the bialgebra case. These do not readily generalize to arbitrary bialgebroids with balancing subalgebra. Thus we consider usual 2-cocycles for bialgebroids, but consider the effect of twisting on the balancing subalgebra.

\begin{definition}
  Let $H$ be a left associative $A$-bialgebroid
  with balancing subalgebra $\mathcal{B}\subset H\otimes_\genfd H$
  such that $\pi(\mathcal{B})\subset H\times_A H$.
  An element $\mathcal{F}\in H\times_A H$ 
  is called a 2-{\bf cocycle} if the equation 
\begin{equation}\label{eq:2coc}
[(\Delta\otimes_A \mathrm{id})(\mathcal{F})](\mathcal{F}\otimes_\genfd 1)
=
[(\mathrm{id}\otimes_A\Delta)(\mathcal{F})](1\otimes_\genfd\mathcal{F})  
\end{equation}
holds in $H\otimes_A H\otimes_A H$.
2-cocycle $\mathcal{F}$ is {\bf counital} if $(\id\otimes_A\epsilon)\mathcal{F} = 1 = (\epsilon\otimes_A\id)\mathcal{F}$. If we write $\mathcal{F} =
\sum_i F^{1i}\otimes F^{2i}:=F^1\otimes F^2$, then the counitality can be rewritten as $\beta(\epsilon(F^2))F^1 = 1 = \alpha(\epsilon(F^1))F^2$.
\end{definition}
This equation~(\ref{eq:2coc}) makes sense by $\mathcal{F}\in H\times_A H$. The case $\pi(\mathcal{B})\subset H\times_A H$ is by Proposition~\ref{prop:equivifBexists} not quite a novel case of a bialgebroid. Still, we are now interested in a recipe for the change of a concrete balancing subalgebra under twisting. Following Xu, for $a\in A$ we define
\begin{equation}\label{eq:alphaFbetaF}
  \beta_{\mathcal{F}}(a) := \beta(F^2\blacktriangleright a)F^1,\,\,\,\,\,\,\,
  \alpha_{\mathcal{F}}(a) := \alpha(F^1\blacktriangleright a)F^2.
\end{equation}
Xu has proved~\cite{xu} that the twisted product $\star_{\mathcal{F}}$ on $A$ defined by $a\star_{\mathcal{F}} b = \alpha_{\mathcal{F}}(a)\beta_{\mathcal{F}}(b)$ is associative and unital. For the $\mathcal{F}$-twisted base algebra $A_{\mathcal{F}} = (A,\star_{\mathcal{F}})$ maps $\alpha_{\mathcal{F}}:A_{\mathcal{F}}\to H$ and $\beta_{\mathcal{F}}:A_{\mathcal{F}}^{\mathrm{op}}\to H$ are morphisms of $\genfd$-algebras with mutually commuting images. In particular, $H$ becomes an $A_{\mathcal{F}}$-bimodule and an $A_{\mathcal{F}}\otimes A_{\mathcal{F}}^{\mathrm{op}}$-ring; use $H^{\mathcal{F}}$ to emphasize the twisted structures. Xu has further shown that  
\begin{equation}\label{eq:FIFeqIA}
\mathcal{F}(\beta_{\mathcal{F}}(a)\otimes 1 - 1\otimes\alpha_{\mathcal{F}}(a))
\in I_A.
\end{equation}
Define $I_{\mathcal{F}}$ as the right ideal in $H\otimes_\genfd H$ generated by all elements of the form $\beta_{\mathcal{F}}(a)\otimes 1 - 1\otimes\alpha_{\mathcal{F}}(a)$. Then~(\ref{eq:FIFeqIA}) implies $\mathcal{F}I_{\mathcal{F}}\subset I_A$. One says that $\mathcal{F}$ is {\bf invertible} if there is an element $\tilde{\mathcal{F}}^{-1}\in H\otimes_\genfd H$ such that $\tilde{\mathcal{F}}^{-1} I_A\subset I_{\mathcal{F}}$ and for $\mathcal{F}^{-1}:=\tilde{\mathcal{F}}^{-1}+I_{\mathcal{F}}$ the identities $\mathcal{F}\mathcal{F}^{-1}= 1\otimes_\genfd 1 + I_A$ and $\mathcal{F}^{-1}\mathcal{F} = 1\otimes_\genfd 1 + I_{\mathcal{F}}$ hold. Denote also by $\tilde{\mathcal{F}}\in H\otimes_\genfd H$ any representative of $\mathcal{F}$. This is not the original definition of invertibility, but it is equivalent to it~\cite{twantip}. It follows that $\mathcal{F} I_{\mathcal{F}} = I_A$ (and $\tilde{\mathcal{F}}I_{\mathcal{F}}=I_A$). Clearly, $H^{\mathcal{F}}\otimes_{A_{\mathcal{F}}}H^{\mathcal{F}} = H\otimes_\genfd H/I_{\mathcal{F}}$. 
If we define $\Delta_{\mathcal{F}}(h) := \mathcal{F}^{-1}\Delta(h)\mathcal{F}:H^{\mathcal{F}}\to H^{\mathcal{F}}\otimes_{A_{\mathcal{F}}}H^{\mathcal{F}}$, this map of $A_{\mathcal{F}}$-bimodules is coassociative due the 2-cocycle property, with counit $\epsilon_{\mathcal{F}}=\epsilon$. Notice that $\Delta_{\mathcal{F}}(h)I_{\mathcal{F}} = \mathcal{F}^{-1}\Delta(h)\mathcal{F}I_{\mathcal{F}}\subset\mathcal{F}^{-1}\Delta(h)I_A \subset\mathcal{F}^{-1}I_A = I_{\mathcal{F}}$. In other words, $\mathrm{Im}\Delta_{\mathcal{F}}\subset H^{\mathcal{F}}\times_{A_{\mathcal{F}}}H^{\mathcal{F}}$, where the factorwise multiplication is well defined. Conjugation with $\mathcal{F}$ (in the sense up to corresponding ideals) can easily be checked to preserve the multiplicativity property of $\Delta_{\mathcal{F}}$. Thus Xu obtains a new twisted $A_{\mathcal{F}}$-bialgebroid $H^{\mathcal{F}}$ from the old $A$-bialgebroid $H$. We want to modify this a bit to allow for a balancing subalgebra.
For this we first describe twisted Takeuchi product in terms
of the original one.
Suppose $\sum_i b_i \otimes_A b'_i\in H\times_A H$, that is 
 $\sum_i \left(b_i\beta_{\mathcal{F}}(a)\otimes_\genfd b'_i -
b_i\otimes_A b'_i\alpha(a) \right) \in I_A$. Then 
$$\begin{array}{lcl}
    \sum_i b_i F^1\beta_{\mathcal{F}}(a)\otimes_A b'_i F^2 &\stackrel{(\ref{eq:alphaFbetaF})}=&
    \sum_i b_i \beta(F^1_{(2)} F'^2\triangleright a) F^1_{(1)} F'^1\otimes_A b'_i F^2\\&\stackrel{(\ref{eq:2coc})}=&\sum_i b_i \beta(F^2_{(1)} F'^1\triangleright a) F^1\otimes_A b'_i F^2_{(2)}F'^2 \\&=&
    \sum_i b_i F^1\otimes_A b'_i \alpha(F^2_{(1)}F'^1\triangleright a)F^2_{(2)}F'^2
\\ &=&\sum_i b_i F^1\otimes_A b'_i F^2\alpha_{\mathcal{F}}(a)
\end{array}$$
Transformations in this calculations are allowed because elements in Takeuchi product multiply elements in $H\otimes_A H$ from the left, and the maps
$(h\otimes_A g)\mapsto \beta(g\triangleright a) h$ and
$(h\otimes_A g)\mapsto \alpha(h\triangleright a)g$ are well defined. 
We obtain $\sum_i b_i F^1\beta_{\mathcal{F}}(a)\otimes_\genfd b'_i F^2 -
  \sum_i b_i F^1 \otimes_\genfd b'_i F^2\alpha_{\mathcal{F}}(a)\in I_A$.
Multiplying this by $\mathcal{F}^{-1}$ from the left, we obtain
\begin{equation}
  \mathcal{F}^{-1}\left(\sum_i b_i\otimes_A b'_i\right)\mathcal{F}\in
  \mathcal{F}^{-1}I_A = I_{\mathcal{F}}. 
\end{equation}
Therefore, $\mathcal{F}^{-1}(H\times_A H)\mathcal{F}\subset H^{\mathcal{F}}\times_{A_{\mathcal{F}}}H^{\mathcal{F}}$,
and similarly for the converse inclusion, obtaining
$$
\mathcal{F}^{-1}(H\times_A H)\mathcal{F} =
H^{\mathcal{F}}\times_{A_{\mathcal{F}}}H^{\mathcal{F}}
$$
Since $I_A$ is right ideal and $\mathcal{F}\mathcal{F}^{-1} = 1\otimes_\genfd 1+I_A$, there are inclusions $I_A\mathcal{F}^{-1} \subset I_A = I_A\mathcal{F}\mathcal{F}^{-1}\subset I_A\mathcal{F}^{-1}$, hence $I_A\mathcal{F}^{-1}=I_A$. Denote the new projection $\pi_{\mathcal{F}}:H\otimes_\genfd H\to (H\otimes_\genfd H)/I_{\mathcal{F}}$. Then define the twisted balancing subalgebra by $\mathcal{B}_{\mathcal{F}} := \pi^{-1}(\mathcal{F}^{-1}\pi(\mathcal{B})\mathcal{F})$. Then $I_{\mathcal{F}}\cap\mathcal{B}_{\mathcal{F}} = (\mathcal{F}^{-1}I_A)\cap\mathcal{B}_{\mathcal{F}} = \mathcal{F}^{-1}(I_A\cap\mathcal{B})\mathcal{F}$. It is a conjugate of a two-sided ideal within algebra $H\otimes_\genfd H$, hence itself a two-sided ideal. If $\Delta\colon H\to\mathcal{B}/(I_A\cap\mathcal{B})$ is a morphism of $\genfd$-algebras, then clearly $\mathcal{F}^{-1}\Delta(-)\mathcal{F}\colon H^{\mathcal{F}}\to\mathcal{B}_{\mathcal{F}}/(I_{\mathcal{F}}\cap\mathcal{B}_{\mathcal{F}})$ is. Thus we obtain a twisted bialgebroid with a balancing subalgebra.

Xu~\cite{xu} does not consider the antipode. A general proof that Drinfeld-Xu twist can be used to twist the antipode by using a canonical formula has been missing for 20 years due technical difficulties resolved only in~\cite{twantip}, for Hopf algebroids with an invertible antipode in the sense of B\"ohm and Szlach\'anyi~\cite{BohmSzlach}. We leave the extension of twisting to antipode in the setting of balancing subalgebras to a future treatment.

{\bf Acknowledgements.} Z. \v{S}koda has been partly supported by grant GA~18-00496S of the Czech Science Foundation. Authors warmly thank Gabriella B\"ohm for reading an early version of the article and suggesting, along with the crucial details, the class of examples of Hopf algebroids with balancing subalgebra constructed from weak Hopf algebras. 

\begin{footnotesize}

\end{footnotesize}

\begin{thebibliography}{99}
\bibitem{bohmnew}
{\sc G. B\"{o}hm}, {\em An alternative notion of Hopf algebroid},
in: Hopf algebras in noncommutative geometry and physics, 
Lecture Notes in Pure and Appl. Math. 239 (2005) 31--53,
\href{https://arxiv.org/abs/math/0311244}{arXiv:math.QA/0311244}.
\bibitem{bohmHbk}
{\sc G. B\"ohm}, {\em Hopf algebroids}, in Handbook of Algebra, 
Vol. 6, ed. by M.~Hazewinkel, Elsevier 2009, 173--236, \href{https://arxiv.org/abs/0805.3806}{arXiv:0805.3806}.
\bibitem{bohmInternal}
{\sc G. B\"ohm}, {\em Internal bialgebroids, entwining structures 
and corings}, AMS Contemp. Math. 376 (2005) 207--226, \href{https://arxiv.org/abs/math/0311244}{arXiv:math.QA/0311244}.
\bibitem{bohmCaenJans2011}
{\sc G. B\"ohm, S. Caenepeel, K. Janssen}, {\em Weak bialgebras and monoidal categories}, Comm. Alg. 39:12 (2011) 4584--4607.
\bibitem{bohmnillszlachwHa}
{\sc G. B\"ohm, F. Nill, K. Szlach\'anyi}, {\em Weak Hopf algebras. I. Integral theory and $C^\ast$-structure}, J. Algebra 221 (1999), no. 2, 385--438, \href{https://arxiv.org/abs/math/9805116}{arXiv:math/9805116}. 
\bibitem{BohmSzlach}
{\sc G. B\"ohm, K. Szlach\'anyi}, {\em Hopf algebroids with
  bijective antipodes: axioms, integrals and duals}, Comm. Alg. 32 (11) (2004) 4433--4464, \href{https://arxiv.org/abs/math/0305136}{arXiv:math.QA/0305136}.
\bibitem{BrzMilitaru}
{\sc T. Brzezi\'nski, G. Militaru}, {\em Bialgebroids,
$\times_A$-bialgebras and duality},  J. Alg. 251: 279--294 (2002)
\href{https://arxiv.org/abs/math/0012164}{arXiv:math.QA/0012164}.
\bibitem{BrzWis}
{\sc T. Brzezi\'nski, R. Wisbauer}, Corings and comodules. 
\emph{London Math. Soc. Lect. Note Ser.,} 309, Cambridge Univ. Press 2003.
\bibitem{daystreetHopfalgd}
{\sc B. Day, R. Street}, {\em Monoidal bicategories and Hopf algebroids}, Advances in Mathematics 129, 1 (1997) 99--157.
\bibitem{doninmudrov}
{\sc J. Donin, A. Mudrov}, {Quantum groupoids and dynamical categories}, 
J. Algebra~\textbf{296}, 348 (2006) \href{https://arxiv.org/abs/math/0311316}{arXiv:math.QA/0311316}.
\bibitem{daystreetqcat}
{\sc R. Street, B. Day}, {\em Quantum categories, star autonomy, and quantum groupoids}, in “Galois Theory, Hopf Algebras, and Semiabelian Categories”, Fields Institute Communications 43 (American Math. Soc. 2004) 187--226; \href{https://arxiv.org/abs/math/0301209}{arXiv:math.CT/0301209}.
\bibitem{ldWeyl}
{\sc N. Durov, S. Meljanac, A. Samsarov, Z. \v{S}koda}, {\em A
universal formula for representing Lie algebra generators as
formal power series with coefficients in the Weyl algebra},
J. Alg. 309, n. 1, 318--359 (2007) \href{https://arxiv.org/abs/math/0604096}{arXiv:math.RT/0604096}.
\bibitem{hovey} 
{\sc M. Hovey}, {\em Morita theory for Hopf algebroids and 
  presheaves of groupoids},  Amer. J. Math. 124:6, 1289--1318 (2002)
\bibitem{kadszl}
{\sc L. Kadison, K. Szlach\'anyi}, Bialgebroid actions on depth two extensions and duality, Adv. Math. 179, Issue 1 (2003) 75--121.
\bibitem{kow}
{\sc N. Kowalzig}, {\em Batalin-Vilkovisky algebra structures on (Co)Tor and Poisson bialgebroids},  J. Pure Appl. Alg. 219:9 (2015) 3781--3822, \href{https://arxiv.org/abs/1305.2992}{arXiv:1305.2992}.
\bibitem{tajronkov}
{\sc T. Juri\'c, D. Kova\v{c}evi\'c, S. Meljanac},  
{\em $\kappa$-deformed phase space, Hopf algebroid and twisting}, 
SIGMA \textbf{10}, 106 (2014), \href{https://arxiv.org/abs/1402.0397}{arXiv:1402.0397}.
\bibitem{tajron}
{\sc T. Juri\'c, S. Meljanac, R. \v{S}trajn},
{\em $\kappa$-Poincar\'e-Hopf algebra and 
Hopf algebroid structure of phase space from twist}, 
Phys. Lett. A 377 (2013), pp. 2472-2476, \href{https://arxiv.org/abs/1303.0994}{arXiv:1303.0994}; 
{\em Twists, realizations and Hopf algebroid structure 
of kappa-deformed phase space}, Int. J. Mod. Phys. A 29 (2014) 1450022, 
\href{https://arxiv.org/abs/1305.3088}{arXiv:1305.3088}.
\bibitem{lu}
{\sc J-H. Lu}, {\em Hopf algebroids and quantum groupoids}, 
Int. J. Math. 7 (1996) 47--70, \href{https://arxiv.org/abs/math/9505024}{arXiv:q-alg/9505024}.
\bibitem{kappaPLB}
{\sc J. Lukierski, Z. \v{S}koda, M. Woronowicz}, {\em Deformed covariant quantum phase spaces as Hopf algebroids}, Phys. Lett. B 750, 401-406 (2015) \href{https://arxiv.org/abs/1507.02612}{arXiv:1507.02612}.
\bibitem{JLZSMWa}
{\sc J. Lukierski, Z. \v{S}koda, M. Woronowicz}, {\em On Hopf algebroid structure of kappa-deformed Heisenberg algebra}, Physics of Atomic Nuclei 80:3 (2017) 569--578, 
\href{https://arxiv.org/abs/1601.01590}{arXiv:1601.01590}.
\bibitem{maltsiniotis}
{\sc G. Maltsiniotis}, {\em Groupo\"ides quantiques},
C. R. Acad. Sci. Paris S\'er. I Math. 314 (1992), no. 4, 249
\bibitem{halgoid}
{\sc S. Meljanac, Z. \v{S}koda, M. Stoji\'c}, {\em
Lie algebra type noncommutative phase spaces are Hopf algebroids}, 
Lett. Math. Phys. 107:3, 475--503 (2017), \href{https://arxiv.org/abs/1409.8188}{arXiv:1409.8188}.
\bibitem{montg}
{\sc S. Montgomery}, {\em Hopf algebras and their actions on
rings}, CBMS Regional Conference Series in Mathematics {\bf 82},
AMS 1993.
\bibitem{Ravenel}
{\sc C. Ravenel}, {\em Complex cobordism and the stable homotopy
groups of spheres}, Acad. Press 1986; 2nd ed. AMS 2003.
\bibitem{schauenburg}
{\sc P. Schauenburg}, {\em Bialgebras over noncommutative rings 
and a structure theorem for Hopf bimodules}, Appl.
Categ. Structures 6 (1998), 193--222.
\bibitem{heisd}
{\sc Z. \v{S}koda}, {\em Heisenberg double versus
deformed derivatives},  Int. J. Mod. Phys. A 26, Nos. 27 \& 28 (2011)
4845--4854, \href{https://arxiv.org/abs/0806.0978}{arXiv:0806.0978}.
\bibitem{twantip}
  {\sc Z. \v{S}koda}, {\em Antipode for Drinfeld-Xu twists of Hopf algebroids}, manuscript
\bibitem{stojicPhD}
{\sc M. Stoji\'c}, {\em Completed Hopf algebroids} (in Croatian, {\em Upotpunjeni Hopfovi algebroidi}), Ph. D. Thesis, Zagreb 2017.
\bibitem{stojicindproVect}
  {\sc M. Stoji\'c}, {\em The symmetric monoidal category indproVect of filtered-cofiltered vector spaces}, in preparation
\bibitem{sweedler}
{\sc M. E. Sweedler}, {\em Groups of simple algebras}, Publ. IH\'ES 44
(1975), 79--189.
\bibitem{takeuchi}
{\sc M. Takeuchi}, {\em Groups of algebras over
$A\otimes\bar{A}$}, J. Math. Soc. Japan 29:3 (1977), 459--492.
\bibitem{xu}
{\sc P. Xu}, {\em Quantum groupoids}, 
Commun. Math. Phys., 216:539--581 (2001) 
\href{https://arxiv.org/abs/math/9905192}{arXiv:q-alg/9905192}.

\end{thebibliography}
\end{document}